\documentclass[12pt,leqno]{amsart}
\newtheorem{thm}{Theorem}[section]
\newtheorem{defi}{Definition}[section]
\newtheorem{lem}{Lemma}[section]
\newtheorem{cor}{Corollary}[section]

\newcommand{\be}{\begin{equation}}
\newcommand{\ee}{\end{equation}}
\numberwithin{equation}{section}
\newcommand{\bea}{\begin{eqnarray}}
\newcommand{\eea}{\end{eqnarray}}
\newcommand{\beb}{\begin{eqnarray*}}
\newcommand{\eeb}{\end{eqnarray*}}
\usepackage{amssymb,amsfonts,amsthm,setspace,indentfirst}
\usepackage{float}
\usepackage[dvips]{graphics}
\usepackage{epsfig}
\usepackage{mathrsfs}
\usepackage{rotating,pdflscape}
\begin{document}
\title[On equivalency  of various geometric structures]{\bf{On equivalency  of various geometric structures}}
\author[Absos Ali Shaikh and Haradhan Kundu]{Absos Ali Shaikh$^*$ and Haradhan Kundu}
\date{}
\address{\noindent\newline Department of Mathematics,\newline University of 
Burdwan, Golapbag,\newline Burdwan-713104,\newline West Bengal, India}
\email{aask2003@yahoo.co.in, aashaikh@math.buruniv.ac.in}
\email{kundu.haradhan@gmail.com}
\begin{abstract}
In the literature we see that after introducing a geometric structure by imposing some restrictions on Riemann-Christoffel curvature tensor, the same type structures given by imposing same restriction on other curvature tensors being studied. The main object of the present paper is to study the equivalency of various geometric structures obtained by same restriction imposing on different curvature tensors. In this purpose we present a tensor by combining Riemann-Christoffel curvature tensor, Ricci tensor, the metric tensor and scalar curvature which describe various curvature tensors as its particular cases. Then with the help of this generalized tensor and using algebraic classification we prove the equivalency of different geometric structures (see, Theorem \ref{cl-1} - \ref{thm5.6}, Table \ref{stc-1} and Table \ref{stc-2}).\\
\end{abstract}
%
\maketitle
\noindent{\bf Mathematics Subject Classiﬁcation (2010).} 53C15, 53C21, 53C25, 53C35.
\section{\bf{Introduction}}\label{intro}
%
Let $M$ be a semi-Riemannian manifold of dimension $n\ge 3$, endowed with the semi-Riemannian metric $g$ with signature $(p, n-p)$, $0\le p \le n$. If (i) $p=0$ or $p=n$; (ii) $p=1$ or $p=n-1$, then $M$ is said to be a (i) Riemannian; (ii) Lorentzian manifold respectively. Let $\nabla$, $R$, $S$ and $r$ be the Levi-Civita connection, Riemannian-Christoffel curvature tensor, Ricci tensor and scalar curvature of $M$ respectively. All the manifolds considered here are assumed to be smooth and connected. We note that any two 1-dimensional semi-Riemannian manifolds are locally isometric, and an 1-dimensional semi-Riemannian manifold is a void field. Also for $n=2$, the notions of above three curvatures are equivalent. Hence throughout the study we will confined ourselves with a semi-Riemannian manifold $M$ of dimension $n\ge 3$. In the study of differential geometry there are various theme of research to derive the geometric properties of a semi-Riemannian manifold. Among others ``symmetry'' plays an important role in the study of differential geometry of a semi-Riemannian manifold.\\
\indent As a generalization of manifold of constant curvature, the notion of local symmetry was introduced by Cartan \cite{ca} with a full classification for the Riemann case. A full classification of such notion was given by Cahen and Parker (\cite{CAH}, \cite{CAH1}) for indefinite case. The manifold $M$ is said to be locally symmetric if its local geodesic symmetries are isometry and $M$ is said to be globally symmetric if its geodesic symmetries are extendible to the whole of $M$. Every locally symmetric manifold is globally symmetric but not conversely. For instance, every compact Riemann surface of genus$> 1$ endowed with its usual metric of constant curvature $(-1)$ is locally symmetric but not globally symmetric. We note that the famous Cartan-Ambrose-Hicks theorem implies that $M$ is locally symmetric if and only if $\nabla R = 0$, and any simply connected complete locally symmetric manifold is globally symmetric. 
During the last eight decades the notion of local symmetry has been weakened by many authors in different directions by imposing some restriction(s) on the curvature tensors and introduced various geometric structures, such as recurrency, semisymmetry, pseudosymmetry etc. In differential geometry there are various curvature tensors arise as an invariant of different transformations, e.g., projective (P), conformal (C), concircular (W), conharmonic (K) curvature tensors etc. All these above restrictions are studied by many geometers on various curvature tensors with their classification, existence and applications.\\
\indent In the literature there are many papers where the same curvature restriction is studied with other curvature tensors which are either meaningless or redundant due to their equivalency. Cartan \cite{ca} first studied the local symmetry. In 1958 So\'os \cite{soos} and in 1964 Gupta \cite{BG} studied the symmetry condition on projective curvature tensor, and then in 1967 Reynolds and Thompson \cite{RT} proved that the notions of local symmetry and projective symmetry are equivalent. Also in \cite{DA2} Desai and Amur studied concircular and projective symmetry and showed that these notions are equivalent to Cartan's local symmetry. Again, the study of recurrent manifold was initiated by Ruse (\cite{Ru1}, \cite{Ru2}, \cite{Ru3}) as the Kappa space and latter named as recurrent space by Walker \cite{Ag}. On the analogy, various authors such as Garai \cite{RG}, Desai and Amur \cite{DA1}, Rahaman and Lal \cite{RL} etc. studied the recurrent notion on the projective curvature tensor as well as concircular curvature tensor. However, all these notions are equivalent to the recurrent manifold (\cite{Glodek}, \cite{mik}, \cite{mik2}, \cite{ol}). Recently Singh \cite{singh} studied the recurrent condition on the $\mathcal{M}$-projective curvature tensor, but from our paper (see, Section 6) it follows that such notion is equivalent to the notion of recurrent manifold. The main object of this paper is to prove the equivalency of various geometric structures obtained by the same curvature restriction on different curvature tensors. For this purpose we present a (0,4) tensor $B$ by the linear combination of the Riemann-Christoffel curvature tensor, Ricci tensor, metric tensor and scalar curvature such that the tensor $B$ describes various curvature tensors as its particular cases. The tensors of the form like $B$ (i.e., particular cases of $B$) are said to be $B$-tensors and the set of all $B$-tensors will be denoted by $\mathscr B$. We classify the set $\mathscr B$ with respect to contraction such that in each class, several geometric structures obtained by the same curvature restriction are equivalent.\\
\indent We are mainly interested on those geometric structures which are obtained by imposing restrictions as some operators on various curvature tensors. We will call these restrictions as ``curvature restriction''. The work on this paper assembled such curvature restrictions and we classify them with respect to their linearity and commutativity with contraction and study the results for each class of restrictions together. On the basis of this study we can say that a specific curvature restriction provides us how many different geometric structures arise due to different curvature tensors.\\
\indent In section 2 we present the tensor $B$ and showed various curvature tensors which are introduced already, are particular cases of it. Section 3 deals with preliminaries. In section 4 we classify the curvature restrictions (these are actually generalized or extended or weaker restrictions of symmetry defined by Cartan) and give the definitions of various geometric structures formed by some curvature restrictions. Section 5 is concerned with basic well known results and some basic properties of the tensors $B$. In section 6 we classify the set $\mathscr B$ and calculate the main results on equivalency of various structures. Finally in last section we make conclusion of the whole work.
\section{\bf{The tensor $B$ and $B$-tensors}}\label{B}
%
Again recall that $M$ is an $n (\ge 3)$-dimensional connected semi-Riemannian manifold equipped with the metric $g$. 
We denote by $\nabla, R, S, r$, the Levi-Civita connection, the Riemann-Christoffel curvature tensor, 
Ricci tensor and scalar curvature of $M$ respectively. We define a $(0,4)$ tensor $B$ given by
\bea\label{tensor}
&&B(X_1,X_2, X_3, X_4) = a_0 R(X_1, X_2, X_3, X_4) + a_1 R(X_1, X_3, X_2, X_4)\\\nonumber
						&&\hspace{0.7in}	+ a_2 S(X_2, X_3) g(X_1, X_4) + a_3 S(X_1, X_3) g(X_2, X_4) + a_4 S(X_1, X_2) g(X_3, X_4)\\\nonumber
						&&\hspace{0.7in}	+ a_5 S(X_1, X_4) g(X_2, X_3) + a_6 S(X_2, X_4) g(X_1, X_3) + a_7 S(X_3, X_4) g(X_1, X_2)\\\nonumber
						&&\hspace{0.7in}	+ r\big[a_8 g(X_1, X_4) g(X_2, X_3) + a_9 g(X_1, X_3) g(X_2, X_4) + a_{10} g(X_1, X_2) g(X_3, X_4)\big],
\eea
where $a_i$'s are scalars on M and $X_1, X_2, X_3, X_4\in \chi(M)$, the Lie algebra of all smooth vector fields on M. Recently Tripathi and Gupta \cite{tri} introduced a similar tensor $\mathcal T$ named as $\mathcal T$-curvature tensor, where two terms are absent. Infact, if $a_1 = 0 = a_{10}$, then the tensor $B$ turns out to be the $\mathcal T$-curvature tensor. Hence the tensor $B$ may be called as \textit{extended $\mathcal T$-curvature tensor} and such name is suggested by M. M. Tripathi (personal communication). However, throughout the paper by the tensor $B$ we shall always mean the extended $\mathcal T$-curvature tensor. We note that for different values of $a_i$'s, as given in the following table (Table \ref{tab-B}), the tensor $B$ reduce to various curvature tensors such as (i) Riemann-Christoffel curvature tensor $R$, (ii) Weyl conformal curvature tensor $C$, (iii) projective curvature tensor $P$, (iv) concircular curvature tensor $W$ \cite{yano1}, (v) conharmonic curvature tensor $K$ \cite{ishi}, (vi) quasi conformal curvature tensor $C^*$ \cite{yano4}, (vii) $C'$ curvature tensor \cite{LC}, (viii) pseudo projective curvature tensor $P^*$ \cite{prasad}, (ix) quasi-concircular curvature tensor $W^*$ \cite{prasad2}, (x) pseudo quasi conformal curvature tensor $\tilde W$ \cite{SJ}, (xi) $\mathcal M$-projective curvature tensor \cite{pokh3}, (xii) $\mathcal W_i$-curvature tensor, $i=1,2,...,9$ (\cite{pokh1}, \cite{pokh2}, \cite{pokh3}), (xiii) $\mathcal W^*_i$-curvature tensor, $i=1,2,...,9$ (\cite{pokh3}) and (xiv) $\mathcal T$-curvature tensor \cite{tri}.
\begin{table}[H]
\begin{center}
{\footnotesize\begin{tabular}{|c|c|c|c|c|c|c|c|c|c|c|c|}\hline
\textbf{Tensor}& $a_0$ & $a_1$ & $a_2$ & $a_3$ & $a_4$ & $a_5$ & $a_6$ & $a_7$ & $a_8$ & $a_9$ & $a_{10}$\\\hline

$R$ &               $1$  & $0$ & $0$ & $0$ & $0$ & $0$ & $0$ & $0$ & $0$ & $0$ & $0$\\\hline

$C$ &               $1$  & $0$ & $-\frac{1}{n-2}$ & $\frac{1}{n-2}$ & $0$ & $-\frac{1}{n-2}$ & $\frac{1}{n-2}$ & $0$ & $\frac{1}{(n-1)(n-2)}$ & $-\frac{1}{(n-1)(n-2)}$ & $0$\\\hline

$P$ &               $1$  & $0$ & $-\frac{1}{n-1}$ & $\frac{1}{n-1}$ & $0$ & $0$ & $0$ & $0$ & $0$ & $0$ & $0$\\\hline

$W$ &        $1$ & $0$ & $0$ & $0$ & $0$ & $0$ & $0$ & $0$ & $-\frac{1}{n(n-1)}$ & $\frac{1}{n(n-1)}$ & $0$\\\hline

$K$ &          $1$  & $0$ & $-\frac{1}{n-2}$ & $\frac{1}{n-2}$ & $0$ & $-\frac{1}{n-2}$ & $\frac{1}{n-2}$ & $0$ & $0$ & $0$ & $0$\\\hline

$C^*$ &             $a_0$  & $0$ & $a_2$ & $-a_2$ & $0$ & $a_2$ & $-a_2$ & $0$ & $-\frac{1}{n}\left(\frac{a_0}{n-1}+2a_2\right)$ & $\frac{1}{n}\left(\frac{a_0}{n-1}+2a_2\right)$ & $0$\\\hline

$C'$ &             $a_0$  & $0$ & $a_2$ & $-a_2$ & $0$ & $a_2$ & $-a_2$ & $0$ & $a_8$ & $-a_8$ & $0$\\\hline

$P^*$ &             $a_0$  & $0$ & $a_2$ & $-a_2$ & $0$ & $0$ & $0$ & $0$ & $-\frac{1}{n}\left(\frac{a_0}{n-1}+a_2\right)$ & $\frac{1}{n}\left(\frac{a_0}{n-1}+a_2\right)$ & $0$\\\hline

$W^*$ &             $a_0$  & $0$ & $0$ & $0$ & $0$ & $0$ & $0$ & $0$ & $\frac{1}{n}\left(\frac{a_0}{n-1}+2 b\right)$ & $\frac{1}{n}\left(\frac{a_0}{n-1}+ 2 b\right)$ & $0$\\\hline

$\widetilde W$ &        $a_0$  & $0$ & $a_2$ & $-a_2$ & $0$ & $a_5$ & $-a_5$ & $0$ & $-\frac{a_0+(n-1)(a_2+a_5)}{n(n-1)}$ & $\frac{a_0+(n-1)(a_2+a_5)}{n(n-1)}$ & $0$\\\hline

$\mathcal M$ &      $1$  & $0$ & $-\frac{1}{2(n-1)}$ & $\frac{1}{2(n-1)}$ & $0$ & $-\frac{1}{2(n-1)}$ & $\frac{1}{2(n-1)}$ & $0$ & $0$ & $0$ & $0$\\\hline

$\mathcal W_0$ &    $1$  & $0$ & $-\frac{1}{n-1}$ & $0$ & $0$ & $0$ & $\frac{1}{n-1}$ & $0$ & $0$ & $0$ & $0$\\\hline

$\mathcal W_0^*$ &  $1$  & $0$ & $\frac{1}{n-1}$ & $0$ & $0$ & $0$ & $-\frac{1}{n-1}$ & $0$ & $0$ & $0$ & $0$\\\hline

$\mathcal W_1$ &    $1$  & $0$ & $-\frac{1}{n-1}$ & $\frac{1}{n-1}$ & $0$ & $0$ & $0$ & $0$ & $0$ & $0$ & $0$\\\hline

$\mathcal W_1^*$ &  $1$  & $0$ & $\frac{1}{n-1}$ & $-\frac{1}{n-1}$ & $0$ & $0$ & $0$ & $0$ & $0$ & $0$ & $0$\\\hline

$\mathcal W_2$ &    $1$  & $0$ & $0$ & $0$ & $0$ & $-\frac{1}{n-1}$ & $\frac{1}{n-1}$ & $0$ & $0$ & $0$ & $0$\\\hline

$\mathcal W_2^*$ &  $1$  & $0$ & $0$ & $0$ & $0$ & $\frac{1}{n-1}$ & $-\frac{1}{n-1}$ & $0$ & $0$ & $0$ & $0$\\\hline

$\mathcal W_3$ &    $1$  & $0$ & $0$ & $-\frac{1}{n-1}$ & $0$ & $\frac{1}{n-1}$ & $0$ & $0$ & $0$ & $0$ & $0$\\\hline

$\mathcal W_3^*$ &  $1$  & $0$ & $0$ & $\frac{1}{n-1}$ & $0$ & $-\frac{1}{n-1}$ & $0$ & $0$ & $0$ & $0$ & $0$\\\hline

$\mathcal W_4$ &    $1$  & $0$ & $0$ & $0$ & $0$ & $0$ & $-\frac{1}{n-1}$ & $\frac{1}{n-1}$ & $0$ & $0$ & $0$\\\hline

$\mathcal W_4^*$ &  $1$  & $0$ & $0$ & $0$ & $0$ & $0$ & $\frac{1}{n-1}$ & $-\frac{1}{n-1}$ & $0$ & $0$ & $0$\\\hline

$\mathcal W_5$ &    $1$  & $0$ & $0$ & $-\frac{1}{n-1}$ & $0$ & $0$ & $\frac{1}{n-1}$ & $0$ & $0$ & $0$ & $0$\\\hline

$\mathcal W_5^*$ &  $1$  & $0$ & $0$ & $\frac{1}{n-1}$ & $0$ & $0$ & $-\frac{1}{n-1}$ & $0$ & $0$ & $0$ & $0$\\\hline

$\mathcal W_6$ &    $1$  & $0$ & $-\frac{1}{n-1}$ & $0$ & $0$ & $0$ & $0$ & $\frac{1}{n-1}$ & $0$ & $0$ & $0$\\\hline

$\mathcal W_6^*$ &  $1$  & $0$ & $\frac{1}{n-1}$ & $0$ & $0$ & $0$ & $0$ & $-\frac{1}{n-1}$ & $0$ & $0$ & $0$\\\hline

$\mathcal W_7$ &    $1$  & $0$ & $-\frac{1}{n-1}$ & $0$ & $0$ & $\frac{1}{n-1}$ & $0$ & $0$ & $0$ & $0$ & $0$\\\hline

$\mathcal W_7^*$ &  $1$  & $0$ & $\frac{1}{n-1}$ & $0$ & $0$ & $-\frac{1}{n-1}$ & $0$ & $0$ & $0$ & $0$ & $0$\\\hline

$\mathcal W_8$ &    $1$  & $0$ & $-\frac{1}{n-1}$ & $0$ & $\frac{1}{n-1}$ & $0$ & $0$ & $0$ & $0$ & $0$ & $0$\\\hline

$\mathcal W_8^*$ &  $1$  & $0$ & $\frac{1}{n-1}$ & $0$ & $-\frac{1}{n-1}$ & $0$ & $0$ & $0$ & $0$ & $0$ & $0$\\\hline

$\mathcal W_9$ &    $1$  & $0$ & $0$ & $0$ & $-\frac{1}{n-1}$ & $\frac{1}{n-1}$ & $0$ & $0$ & $0$ & $0$ & $0$\\\hline

$\mathcal W_9^*$ &  $1$  & $0$ & $0$ & $0$ & $\frac{1}{n-1}$ & $-\frac{1}{n-1}$ & $0$ & $0$ & $0$ & $0$ & $0$\\\hline

$\tau$ & $a_0$ & $0$ & $a_2$ & $a_3$ & $a_4$ & $a_5$ & $a_6$ & $a_7$ & $a_8$ & $a_9$ & $0$\\\hline

\end{tabular}}
\end{center}
\caption{List of $B$-tensors}\label{tab-B}
\end{table}
There may arise some other tensors from the tensor $B$ as its particular cases, which are not introduced so far. We recall that the tensors arising out from the tensor $B$ as its particular cases are called $B$-tensors and the set of all such $B$-tensors will be denoted by $\mathscr B$. It is easy to check that $\mathscr B$ forms a module over $C^{\infty}(M)$, the ring of all smooth functions on $M$. We note that the $B$-tensor pseudo quasi-conformal curvature tensor $\widetilde W$ was studied by Shaikh and Jana in 2006 \cite{SJ}, but the same notion was studied by Prasad et. al. in 2011 \cite{prasad3} as generalized quasi-conformal curvature tensor ($G_{qc}$).
%
\section{\bf{Preliminaries}}\label{pre}
%
Let us now consider a connected semi-Riemannian manifold $M$ of dimension n($\ge 3$). Then for two (0, 2) tensors $A$ and $E$, the 
Kulkarni-Nomizu product (\cite{D0}, \cite{DGHS}, \cite{GLOG}, \cite{GLOG1}, \cite{Kow2}) $A\wedge E$ is given by
\bea\label{KN}
(A \wedge E)(X_1,X_2,X_3,X_4)&=&A(X_1,X_4)E(X_2,X_3) + A(X_2,X_3)E(X_1,X_4)\\\nonumber
&&-A(X_1,X_3)E(X_2,X_4) - A(X_2,X_4)E(X_1,X_3),
\eea
where $X_1, X_2, X_3, X_4\in \chi(M)$.\\
A tensor $D$ of type (1,3) on $M$ is said to be generalized curvature tensor (\cite{D11}, \cite{D2}, \cite{D10}), if
\beb
&(i)&D(X_1,X_2)X_3+D(X_2,X_1)X_3=0,\\
&(ii)&D(X_1,X_2,X_3,X_4)=D(X_3,X_4,X_1,X_2),\\
&(iii)&D(X_1,X_2)X_3+D(X_2,X_3)X_1+D(X_3,X_1)X_2=0,
\eeb
where $D(X_1,X_2,X_3,X_4)=g(D(X_1,X_2)X_3,X_4)$, for all $X_1,X_2,$ $X_3,X_4$. Here we denote the same symbol $D$ for both generalized curvature tensor of type (1,3) and (0,4). Moreover if $D$ satisfies the second Bianchi identity i.e.,
$$(\nabla_{X_1}D)(X_2,X_3)X_4+(\nabla_{X_2}D)(X_3,X_1)X_4+(\nabla_{X_3}D)(X_1,X_2)X_4=0,$$
then $D$ is called a proper generalized curvature tensor. We note that a linear combination of generalized curvature tensors over $C^{\infty}(M)$ is again a generalized curvature tensor but it is not true for proper generalized curvature tensors, in general. However, if the linear combination is taken over $\mathbb R$, then it is true.\\
%
\indent Now for any (1,3) tensor $D$ (not necessarily generalized curvature tensor) and given two vector fields $X,Y\in\chi(M)$, one can define an endomorphism $\mathscr{D}(X,Y)$ by
$$\mathscr{D}(X,Y)(Z)=D(X,Y)Z, \ \ \forall\mbox{$Z\in\chi(M)$}.$$
Again, if $X,Y\in\chi(M)$ then for a (0,2) tensor $A$, one can define two endomorphisms $\mathscr{A}$ and $X \wedge_A Y$, by (\cite{D11}, \cite{D2}, \cite{D10})
$$g(\mathscr{A}(X),Y)=A(X,Y),$$
$$(X \wedge_A Y)Z = A(Y,Z)X - A(X,Z)Y, \ \forall\ \mbox{$Z \in \chi(M)$}.$$
Now for a $(0,k)$-tensor $T$, $k\geq 1$, and an endomorphism $\mathscr H$, one can operate $\mathscr H$ on $T$ to produce the 
tensor $\mathscr H T$, given by (\cite{D11}, \cite{D2}, \cite{D10})
\beb\label{rdot}
(\mathscr{H} T)(X_1,X_2,\cdots,X_k) = -T(\mathscr{H}X_1,X_2,\cdots,X_k) - \cdots - T(X_1,X_2,\cdots,\mathscr{H}X_k).
\eeb
\indent We consider that the operation of $\mathscr H$ on a scalar is zero. In particular, $\mathscr H$ may be $\mathscr{D}(X,Y)$, $X \wedge_A Y$, $\mathscr{A}$ etc. In particular for $\mathscr H = \mathscr{D}(X,Y)$ and $\mathscr{H} = (X \wedge_A Y)$, we have (\cite{D11}, \cite{D2}, \cite{D10}, \cite{tac})
\beb\label{ddt}
(\mathscr D(X,Y) T)(X_1,X_2,\cdots,X_k)=-T(\mathscr D(X,Y)(X_1),X_2,\cdots,X_k) - \cdots - T(X_1,X_2,\cdots,\mathscr D(X,Y)(X_k))\\\nonumber
=- T(D(X,Y)X_1,X_2,\cdots,X_k) - \cdots - T(X_1,X_2,\cdots, D(X,Y)X_k),
\eeb
\beb\label{qgr}
((X \wedge_A Y) T)(X_1,X_2,\cdots,X_k)= -T((X \wedge_A Y)X_1,X_2,\cdots,X_k) - \cdots - T(X_1,X_2,\cdots,(X \wedge_A Y)X_k)\\\nonumber
= A(X, X_1) T(Y,X_2,\cdots,X_k) + \cdots + A(X, X_k) T(X_1,X_2,\cdots,Y)\\\nonumber
- A(Y, X_1) T(X,X_2,\cdots,X_k) - \cdots - A(Y, X_k) T(X_1,X_2,\cdots,X),
\eeb
where $X, Y, X_i \in \chi(M)$, $i = 1,2,\cdots,k$.\\
We denote the above tensor $(\mathscr D(X,Y) T)(X_1,X_2,\cdots,X_k)$ as $D\cdot T(X_1,X_2,\cdots,X_k,X,Y)$ and the tensor $((X \wedge_A Y) T)(X_1,X_2,\cdots,X_k)$ as $Q(A,T)(X_1,X_2,\cdots,X_k,X,Y)$.\\
\indent For an 1-form $\Pi$ and a vector field $X$ on $M$, we can define an endomorphism $\Pi_{_X}$ as
$$\Pi_{_X}(X_1) = \Pi(X_1)X, \ \mbox{$\forall X_1\in \chi(M)$}.$$
Then we can define the tensor $\Pi_{_X} T$ as follows:
\beb\label{pidot}
&&(\Pi_{_X} T)(X_1,X_2, \cdots, X_k)\\\nonumber
&&= -T(\Pi_{_X}(X_1),X_2, \cdots, X_k) - \cdots - T(X_1,X_2, \cdots, \Pi_{_X}(X_k)),\\\nonumber
&&= -\Pi(X_1)T(X,X_2, \cdots, X_k) -\Pi(X_2)T(X_1,X, \cdots, X_k)- \cdots -\Pi(X_k)T(X_1,X_2, \cdots, X),
\eeb
$\forall X, X_i\in\chi(M), i= 1, 2, \cdots,k$.
%
\section{\bf{Some geometric structures defined by curvature related operators}}\label{structures}
In this section we discuss some geometric structures which arise by some curvature restrictions on a semi-Riemannian manifold. We are mainly interested on those geometric structures which are obtained by some curvature restrictions imposed on $B$-tensors by means of some operators, e.g., symmetry, recurrency, pseudosymmetry etc. These operators are linear over $\mathbb{R}$ and may or may not be linear over $C^{\infty}(M)$ and thus called as $\mathbb{R}$-linear operators or simply linear operators. The linear operators which are not linear over $C^{\infty}(M)$, said to be \emph{operators of the 1st type} and which are linear over $C^{\infty}(M)$, said to be \emph{operators of the 2nd type}. Some important 1st type operators are symmetry, recurrency, weakly symmetry (in the sense of Tam\'assy and Binh) etc. and some important 2nd type operators are semisymmetry, Deszcz pseudosymmetry, Ricci generalized pseudosymmetry etc. We denote the set of all tensor fields on $M$ of type $(r,s)$ by $\mathcal T^k_s$ and we take $\mathcal L$ as any $\mathbb R$-linear operator such that the operation of $\mathcal L$ on $T \in \mathcal T^k_s$ is denoted by $\mathcal L\ T$.\\
\indent Another classification of such linear operators may be given with respect to their extendibility. Actually these operators are imposed on $(0,4)$ curvature tensors but the defining condition of some of them can not be extended to any $(0,k)$ tensor, e.g., symmetry, semisymmetry, weak symmetry (all three types) operators are extendible but weakly generalized recurrency, hyper generalized recurrency operators are not extendible. Again extendible operators are classified into two subclasses, (i) operators commute with contraction or commutative and (ii) operators not commute with contraction or non-commutative., e.g., symmetry, semisymmetry operators are commutative but weak symmetry operators are non-commutative. Throughout this paper by a commutative or non-commutative operator we mean a linear operator which commutes or not commutes with contraction. The tree diagram of the classification of linear operators imposed on (0,4) curvature tensors is given by:\\

\vspace{0.1in}
\noindent $\put(20,0){\framebox(120,30){operators of 1st type}}       \put(300,0){\framebox(120,30){operators of 2nd type}}$
\noindent $\put(80,-16){\vector(0,1){15}}    \put(380,-16){\vector(0,1){15}}
\put(80,-16){\line(1,0){300}}
\put(225,-31){\vector(0,1){15}}$\\
$\put(120,0){\framebox(200,50){$\begin{array}{c}$\textbf{Class of linear operators defined}$\\$\textbf{on (0,4) curvature tensors}$\end{array}$}}
\put(225,-1){\vector(0,-1){15}}
\put(80,-16){\line(1,0){300}}
\put(80,-16){\vector(0,-1){15}}    \put(380,-16){\vector(0,-1){15}}
\put(20,-61){\framebox(182,30){non-extendible for any (0,k) tensor}}
\put(290,-61){\framebox(160,30){extendible for any (0,k) tensor}}
\put(340,-62){\vector(0,-1){15}}
\put(160,-77){\line(1,0){240}}
\put(160,-77){\vector(0,-1){15}}    \put(400,-77){\vector(0,-1){15}}
\put(110,-132){\framebox(140,40){$\begin{array}{c}$commute with contraction$\\$or commutative$\end{array}$}}
\put(320,-132){\framebox(160,40){$\begin{array}{c}$not commute with contraction$\\$or non-commutative$\end{array}$}}
$
\vspace{0.1in}
\begin{defi}\cite{ca} 
Consider the covariant derivative operator $\nabla_X : \mathcal T^0_k \rightarrow \mathcal T^0_{k+1}$. A semi-Riemannian manifold is said to be \textit{$T$-symmetric} if $\nabla_X T = 0$, for all $X\in \chi(M)$, $T\in \mathcal T^0_{k}$.
\end{defi}
\vspace{-0.3cm}
\indent Obviously this operator is of 1st type and commutative. The condition for $T$-symmetry is written as $\nabla T = 0$.
\begin{defi}$($\cite{Ag}, \cite{Ru1}, \cite{Ru2}, \cite{Ru3}$)$ 
Consider the operator $\kappa_{(X,\Pi)}: \mathcal T^0_k \rightarrow \mathcal T^0_{k+1}$ defined by $\kappa_{(X,\Pi)} T = \nabla_X T - \Pi(X)\otimes T$, $\otimes$ is the tensor product, $\Pi$ is an 1-form and $T\in \mathcal T^0_k$. A semi-Riemannian manifold is said to be \textit{$T$-recurrent} if $\kappa_{(X,\Pi)}T = 0$ for all $X\in \chi(M)$ and some 1-form $\Pi$, called the associated 1-form or the 1-form of recurrency.
\end{defi}
\vspace{-0.2cm}
\indent Obviously this operator is of 1st type and commutative. The condition for $T$-recurrency is written as $\nabla T - \Pi\otimes T = 0$ or simply $\kappa T = 0$.\\
\indent Keeping the commutativity property, we state some generalization of symmetry operator and recurrency operator which are respectively said to be \textit{symmetric type} operator and \textit{recurrent type} operator. For this purpose we denote the $s$-th covariant derivative as
$$\nabla_{X_1}\nabla_{X_2}\cdots\nabla_{X_s} = \nabla^{s}_{X_1 X_2\cdots X_s}.$$
Now the operator
$$L^s_{X_1 X_2 \cdots X_s} = \sum_{\sigma} \alpha_{\sigma}\nabla^{s}_{X_{\sigma(1)} X_{\sigma(2)}\cdots X_{\sigma(s)}}$$
is called a \textit{symmetric type operator of order $s$}, where $\sigma$ is permutation over $\{1,2,...,s\}$ and the sum is taken over the set of all permutations over $\{1,2,...,s\}$ and $\alpha_{\sigma}$'s are some scalars not all together zero.\\
\indent A manifold is called \textit{$T$-symmetric type of order $s$} if
\be\label{eq4.4}
L^s_{X_1 X_2 \cdots X_s} T = 0 \ \ \ \forall\  X_1, X_2, \cdots X_s \in\chi(M) \mbox{ and some scalars } \alpha_{\sigma}.
\ee
The scalars $\alpha_{\sigma}$'s are called the associated scalars. We denote the condition for $T$-symmetry type of order $s$ is written as $L^s T = 0$.\\
Again for some $(0,i)$ tensors $\Pi^i_{\sigma}$ (not all together zero), $i=0,1,2,\cdots,s$ and all permutations $\sigma$ over $\{1,2,...,s\}$  (i.e., $\Pi^0_{\sigma}$ are scalars), the operator
\beb
\kappa^s_{X_1 X_2 \cdots X_s} &=&\sum_{\sigma}\left[ \Pi^0_{\sigma}  \nabla^{s}_{X_{\sigma(1)} X_{\sigma(2)}\cdots X_{\sigma(s)}}\right.\\
															&& + \Pi^1_{\sigma}(X_{\sigma(1)})  \nabla^{s-1}_{X_{\sigma(2)} X_{\sigma(3)}\cdots X_{\sigma(s)}}\\
															&& + \Pi^2_{\sigma}(X_{\sigma(1)},X_{\sigma(2)})  \nabla^{s-2}_{X_{\sigma(3)} X_{\sigma(4)}\cdots X_{\sigma(s)}}\\
															&& + \cdots \ \ \ \cdots \ \ \ \cdots \ \ \ \cdots\\
															&& + \Pi^{s-1}_{\sigma}(X_{\sigma(1)},X_{\sigma(2)},\cdots, X_{\sigma(s-1)})  \nabla_{X_{\sigma(s)}}\\
															&& \left. + \Pi^s_{\sigma}(X_{\sigma(1)},X_{\sigma(2)},\cdots, X_{\sigma(s)})I_d\right],
\eeb
$I_d$ is the identity operator, is called a \textit{recurrent type operator of order $s$}.\\
\indent A manifold is called \textit{$T$-recurrent type of order $s$} if it satisfies
\be\label{eq4.5}
\kappa^s_{X_1 X_2 \cdots X_s} T = 0, \ \ \forall X_1, X_2, \cdots X_s \in\chi(M) \mbox{ and some $i$-forms $\Pi^i_{\sigma}$'s}, i=0,1,2,\cdots,s.
\ee
The $i$-forms $\pi^{i}_{\sigma}$'s are called the associated $i$-forms. The $T$-recurrency condition of order $s$, will be simply written as $\kappa^s T = 0$.\\
\indent As recurrency is a generalization of symmetry, likewise, the recurrent type condition is a generalization of symmetric type condition. We note that these symmetric type and recurrent type operators are, generally, of 1st type but some of them may be of second type. For example, the semisymmetric operator $(\nabla_X\nabla_Y-\nabla_Y\nabla_X)$  is symmetric type as well as recurrent type and also of 2nd type operator.\\
\indent Another way to generalize recurrency there are some other geometric structures defined as follows:
\begin{defi}\cite{DUB} 
Consider the operator $G\kappa_{(X,\Pi,\Phi)}: \mathcal T^0_4 \rightarrow \mathcal T^0_5$ defined by
$$G\kappa_{(X,\Pi,\Phi)} T = \nabla_X T - \Pi(X)\otimes T - \Phi(X)\otimes G,$$
$\Pi$ and $\Phi$ are 1-forms and $T$ is a $(0,4)$ tensor. A semi-Riemannian manifold is said to be generalized $T$-recurrent if $G\kappa_{(X,\Pi,\Phi)}T = 0$ for all $X\in\chi(M)$ and some 1-forms $\Pi$ and $\Phi$, called the associated 1-forms.
\end{defi}
\vspace{-0.4cm}
\indent Obviously this operator is of 1st type and non-extendible.
\begin{defi}\cite{SP} 
Consider the operator $H\kappa_{(X,\Pi,\Phi)}: \mathcal T^0_4 \rightarrow \mathcal T^0_5$ defined by
$$H\kappa_{(X,\Pi,\Phi)} T = \nabla_X T - \Pi(X)\otimes T - \Phi(X)\otimes g\wedge S,$$
$\Pi$ and $\Phi$ are 1-forms and $T$ is a $(0,4)$ tensor. A semi-Riemannian manifold is said to be hyper-generalized $T$-recurrent if $H\kappa_{(X,\Pi,\Phi)}T = 0$ for all $X\in\chi(M)$ and some 1-forms $\Pi$ and $\Phi$, called the associated 1-forms.
\end{defi}
\indent Obviously this operator is of 1st type and non-extendible.
\begin{defi}\cite{ROY} 
Consider the operator $W\kappa_{(X,\Pi,\Phi)}: \mathcal T^0_4 \rightarrow \mathcal T^0_5$ defined by
$$W\kappa_{(X,\Pi,\Phi)} T = \nabla_X T - \Pi(X)\otimes T - \Phi(X)\otimes S\wedge S,$$
$\Pi$ and $\Phi$ are 1-forms and $T$ is a $(0,4)$ tensor. A semi-Riemannian manifold is said to be weakly generalized $T$-recurrent if $W\kappa_{(X,\Pi,\Phi)}T = 0$ for all $X\in\chi(M)$ and some 1-forms $\Pi$ and $\Phi$, called the associated 1-forms.
\end{defi}
\indent Obviously this operator is of 1st type and non-extendible.
\begin{defi}\cite{ROY1} 
Consider the operator $Q\kappa_{(X,\Pi,\Phi)}: \mathcal T^0_4 \rightarrow \mathcal T^0_5$ defined by
$$Q\kappa_{(X,\Pi,\Phi,\Psi)} T = \nabla_X T - \Pi(X)\otimes T - \Phi(X)\otimes [g\wedge (g + \Psi\otimes \Psi)],$$
$\Pi$, $\Phi$ and $\Psi$ are 1-forms and $T$ is a $(0,4)$ tensor. A semi-Riemannian manifold is said to be quasi generalized $T$-recurrent if $Q\kappa_{(X,\Pi,\Phi, \Psi)}T = 0$ for all $X\in\chi(M)$ and some 1-forms $\Pi$, $\Phi$ and $\Psi$, called the associated 1-forms.
\end{defi}
\indent Obviously this operator is of 1st type and non-extendible.
\begin{defi} 
Consider the operator $S\kappa_{(X,\Pi,\Phi,\Psi,\Theta)}: \mathcal T^0_4 \rightarrow \mathcal T^0_5$ defined by
$$S\kappa_{(X,\Pi,\Phi,\Psi,\Theta)} T = \nabla_X T - \Pi(X)\otimes T - \Phi(X)\otimes G - \Psi(X)\otimes g\wedge S - \Theta(X)\otimes S\wedge S,$$
$\Pi$, $\Phi$, $\Psi$ and $\Theta$ are 1-forms and $T$ is a $(0,4)$ tensor. A semi-Riemannian manifold is said to be super generalized $T$-recurrent if $S\kappa_{(X,\Pi,\Phi,\Psi,\Theta)}T = 0$ for all $X\in\chi(M)$ and some 1-forms $\Pi$, $\Phi$, $\Psi$ and $\Theta$, called the associated 1-forms.
\end{defi}
\indent Obviously this operator is of 1st type and non-extendible.\\
\indent We now state another generalization of local symmetry, given as follows:
\begin{defi}\cite{CHA} 
Consider the operator $CP_{(X,\Pi)}: \mathcal T^0_k \rightarrow \mathcal T^0_{k+1}$ defined by
$$CP_{(X,\Pi)} T = \nabla_X T - 2\Pi(X)\otimes T + \Pi_{_X} T,$$
$\Pi$ is an 1-form and $T\in \mathcal T^0_k$. A semi-Riemannian manifold is said to be Chaki $T$-pseudosymmetric \cite{CHA} if $CP_{(X,\Pi)} T = 0$ for all $X\in\chi(M)$ and some 1-form $\Pi$, called the associated 1-form. 
\end{defi}
\indent Obviously this operator is of 1st type and non-commutative.\\
\indent Again in another way Tam$\acute{\mbox{a}}$ssy and Binh \cite{tb} generalized the recurrent and Chaki pseudosymmetric structures and named it weakly symmetric structure. But there are three types of weak symmetry \cite{SD} which are given below:
\begin{defi}
Consider the operator $W^1_{(X,\stackrel{\sigma}{\Pi})}: \mathcal T^0_k \rightarrow \mathcal T^0_{k+1}$ defined by 
$$W^1_{(X,\stackrel{\sigma}{\Pi})} T = (\nabla_{X}T)(X_2,X_3,...,X_{k+1}) - \sum_{\sigma}\stackrel{\sigma}{\Pi}(X_{\sigma(1)})T(X_{\sigma(2)},X_{\sigma(3)},...,X_{\sigma(k+1)}),$$
$\stackrel{\sigma}{\Pi}$ are 1-forms, $T\in \mathcal T^0_k$ and the sum includes all permutations $\sigma$ over the set $(1,2,...,k+1)$. A semi-Riemannian manifold $M$ is said to be weakly $T$-symmetric of type-I if $W^1_{(X,\stackrel{\sigma}{\Pi})} T = 0,$ for all $X\in\chi(M)$ and some 1-forms $\stackrel{\sigma}{\Pi}$, called the associated 1-forms.
\end{defi}
\indent Obviously this operator is of 1st type and non-commutative.
\begin{defi}
Consider the operator $W^2_{(X,\Phi,\Pi_i)}: \mathcal T^0_k \rightarrow \mathcal T^0_{k+1}$ defined by
\beb
(W^2_{(X,\Phi,\Pi_i)} T)(X_1,X_2,...,X_k) &=& (\nabla_X T)(X_1,X_2,...,X_k)\\\nonumber &-&\Phi(X)T(X_1,X_2,...,X_k)-\sum^k_{i=1}\Pi_i(X_i)T(X_1,X_2,...,\underset{i-th\ place}{X},...,X_k),
\eeb
where $\Phi$ and $\Pi_i$ are 1-forms and $T\in \mathcal T^0_k$. A semi-Riemannian manifold $M$ is said to be weakly $T$-symmetric of type-II if $W^2_{(X,\Phi,\Pi_i)} T = 0,$ for all $X\in\chi(M)$ and some 1-forms $\Phi$ and $\Pi_i$, called the associated 1-forms.
\end{defi}
\indent Obviously this operator is of 1st type and non-commutative.
\begin{defi}
Consider the operator $W^3_{(X,\Phi,\Pi)}: \mathcal T^0_k \rightarrow \mathcal T^0_{k+1}$ defined by 
$$W^3_{(X,\Phi,\Pi)} T = \nabla_X T - \Phi\otimes T  + \pi_{_X}T,$$
where $\Phi$ and $\Pi$ are 1-forms and $T\in \mathcal T^0_k$. A semi-Riemannian manifold $M$ is said to be weakly $T$-symmetric of type-III if $W^3_{(X,\Phi,\Pi)} T = 0,$ for all $X\in\chi(M)$ and two 1-forms $\Phi$ and $\Pi$, called the associated 1-forms.
\end{defi}
\indent Obviously this operator is of 1st type and non-commutative.\\
\indent The weak symmetry of type-II was first introduced by Tam$\acute{\mbox{a}}$ssy and Binh \cite{tb} and the other two types of the weak symmetry can be deduced from the type-II (see, \cite{sti}). Although there is an another notion of weak symmetry introduced by Selberg \cite{sel} which is totally different from this notion and the representation of such a structure by the curvature restriction is unknown till now. However, throughout our paper we will consider the weak symmetry in sense of Tam$\acute{\mbox{a}}$ssy and Binh \cite{tb}.
\begin{defi} 
For a $(0,4)$ tensor $D$ consider the operator $\mathcal D(X,Y): \mathcal T^0_r \rightarrow \mathcal T^0_{r+2}$ defined by
$$(\mathcal D(X,Y) T)(X_1,X_2,\cdots,X_k) = (D\cdot T)(X_1,X_2,\cdots,X_k,X,Y).$$
A semi-Riemannian manifold is said to be $T$-semisymmetric type if $\mathcal D(X, Y) T = 0$ for all $X,Y\in\chi(M)$. This condition is also written as $D\cdot T = 0$.
\end{defi}
\indent Obviously this operator is of 2nd type and commutative or non-commutative according as $D$ is skew-symmetric or not in 3rd and 4th places i.e., $D(X_1,X_2,X_3,X_4) = -D(X_1,X_2,X_4,X_3)$ or not. Especially, if we consider $D = R$, then the manifold is called T-semisymmetric \cite{sz}.
\begin{defi}$($\cite{adm}, \cite{DR1}, \cite{DES}, \cite{des7}$)$ 
A semi-Riemannian manifold is said to be $T$-pseudosymmetric type if $(\sum_i c_i D_i)\cdot T = 0$, where $\sum_i c_i D_i$ is a linear combination of $(0,4)$ curvature tensors $D_i$'s over $C^{\infty}(M)$, $c_i\in C^{\infty}(M)$, called the associated scalars.
\end{defi}
\indent Obviously this operator is of 2nd type and generally commutative or non-commutative according as all $D_i$'s are skew-symmetric or not in 3rd and 4th places. Consider the special cases $(R - L G)\cdot T = 0$ and $(R - L X\wedge_S Y)\cdot T = 0$. These are known as Deszcz T-pseudosymmetric (\cite{adm}, \cite{DR1}, \cite{DES}, \cite{des7}) and Ricci generalized T-pseudosymmetric (\cite{DEF}, \cite{DEF1}) respectively. It is clear that the operator of Deszcz pseudosymmetry is commutative but Ricci generalized pseudosymmetry is non-commutative.
%
\section{\bf{Some basic properties of the tensor $B$}}\label{basic}
%
In this section we discuss some basic well known properties of the tensor $B$.
\begin{lem}\label{lemsk}
An operator $\mathcal L$ is commutative if $\mathcal L g = 0$. Moreover if $\mathcal L$ is an endomorphism then this condition is equivalent to the condition that $\mathcal L$ is skew-symmetric i.e. $g(\mathcal L X, Y) = - g(X, \mathcal L Y)$ for all $X,Y\in\chi(M)$.
\end{lem}
\noindent\textbf{Proof:} If $\mathcal L  g = 0$ then, without loss of generality, we may suppose that T is a (0,2) tensor, and we have
$$\mathcal L(\mathscr{C}(T)) = \mathcal L (g^{ij}T_{ij}) = g^{ij} (\mathcal L T_{ij}) = \mathscr{C}(\mathcal L T),$$
where $\mathscr{C}$ is the contraction operator.
Again if $\mathcal L$ is an endomorphism then for all $X,Y\in\chi(M)$, $\mathcal L g = 0$ implies
$$(\mathcal L g)(X, Y) = -g(\mathcal L X, Y)-g(X, \mathcal L Y) = 0$$
$$\Rightarrow g(\mathcal L X, Y) = - g(X, \mathcal L Y)$$
\begin{center}$\Rightarrow \mathcal L$ is skew-symmetric.\end{center}
From the last part of this lemma we can say that
\begin{lem}\label{lem5.22}
The curvature operator $\mathscr D(X,Y)$ formed by a $(0,4)$ tensor $D$ is commutative if and only if $D$ is skew-symmetric in 3rd and 4th places, i.e., $D(X_1,X_2,X_3,X_4)=-D(X_1,X_2,X_4,X_3)$,  for all $X_1,X_2,X_3,X_4$.
\end{lem}
\begin{lem}\label{lem3.1}
Contraction and covariant derivative operators are commute each other.
\end{lem}
\begin{lem}\label{lem3.5}
$Q(g,T) = G\cdot T$.
\end{lem}
\noindent\textbf{Proof:} For a (0,k) tensor $T$, we have
\beb
Q(g,T)(X_1,X_2,\cdots X_k;X,Y)=((X\wedge_g Y)\cdot T)(X_1,X_2,\cdots X_k).
\eeb
Now $(X\wedge_g Y)(X_1,X_2) = G(X,Y,X_1,X_2)$, so the result follows.
\begin{lem}\label{lem3.6}
Let $D$ be a generalized curvature tensor. Then\\
(1) $D(X_1,X_2,X_1,X_2) = 0$ implies $D(X_1,X_2,X_3,X_4) = 0$,\\
(2) $(\mathcal L D)(X_1,X_2,X_1,X_2) = 0$ implies $(\mathcal L D)(X_1,X_2,X_3,X_4) = 0$, $\mathcal L$ is any linear operator.
\end{lem}
\noindent\textbf{Proof:} The results follows from Lemma 8.9 of \cite{lee} and hence we omit it.\\
%
%
\indent We now consider the tensor $B$ and take contraction on i-th and j-th place and get the $(i$-$j)$-th contraction tensor  $^{^{ij}}S$ for $i, j \in \{1,2,3,4\}$ as
\bea
&&\left\{\begin{array}{l}
^{12}S = (-a_1 + a_2 + a_3 + a_5 + a_6 + n a_7)S + r(a_4 + a_8 + a_9 + n a_{10})g 			= (^{12}p) S + (^{12}q) r g\\
^{13}S = (-a_0 + a_2 + a_4 + a_5 + n a_6 + a_7)S + r(a_3 + a_8 + n a_9 + a_{10})g 			= (^{13}p) S + (^{13}q) r g\\
^{14}S = (a_0 + a_1 + n a_2 + a_3 + a_4 + a_6 + a_7)S + r(a_5 + n a_8 + a_9 + a_{10})g 	= (^{14}p) S + (^{14}q) r g\\
^{23}S = (a_0 + a_1 + a_3 + a_4 + n a_5 + a_6 + a_7)S + r(a_2 + n a_8 + a_9 + a_{10})g 	= (^{23}p) S + (^{23}q) r g\\
^{24}S = (-a_0 + a_2 + n a_3 + a_4 + a_5 + a_7)S + r(a_6 + a_8 + n a_9 + a_{10})g 			= (^{24}p) S + (^{24}q) r g\\
^{34}S = (-a_1 + a_2 + a_3 + n a_4 + a_5 + a_6)S + r(a_7 + a_8 + a_9 + n a_{10})g 			= (^{34}p) S + (^{34}q) r g
\end{array}\right.
\eea
Again contracting all $^{^{ij}}S$ we get $^{ij}r$ for $i, j \in \{1,2,3,4\}$ as
\bea
&&\left\{\begin{array}{l}
^{12}r = {}^{34}r = (-a_1 + a_2 + a_3 + n a_4 + a_5 + a_6+ n a_7+ n a_8+ n a_9 + n^2 a_{10}) r\\
^{13}r = {}^{24}r = (-a_0 + a_2 + n a_3 + a_4 + a_5 + n a_6 + a_7 + n a_8+ n^2 a_9 + n a_{10}) r\\
^{14}r = {}^{23}r = (a_0 + a_1 + n a_2 + a_3 + a_4 + n a_5 + a_6 + a_7 +  n^2 a_8 +  n a_9 +  n a_{10}) r
\end{array}\right.
\eea
\begin{lem}\label{lem3.8}
$(i)$ If $S = 0$, then $B = 0$ if and only if $R = 0$.\\
$(ii)$ If $\mathcal L S = 0$, then $\mathcal L B = 0$ if and only if $\mathcal L R = 0$, where $\mathcal L$ is a commutative 1st type operator and $a_i$'s are constant.\\
$(iii)$ If $\mathcal L S = 0$, then $\mathcal L B = 0$ if and only if $\mathcal L R = 0$, where $\mathcal L$ is a commutative 2nd type operator.
\end{lem}
\begin{lem}\label{lem3.9}
The tensor $B$ is a generalized curvature tensor if and only if
\bea\label{eq5.3}
&&\left\{\begin{array}{c}
a_1 = a_4 = a_7 = a_{10} = 0,\\
a_2 = - a_3 = a_5 = - a_6 \mbox{ and } a_8 = -a_9.
\end{array}\right.
\eea
\end{lem}
\noindent\textbf{Proof:} $B$ is a generalized curvature tensor if and only if
\bea\label{eq3.1}
&&\left\{\begin{array}{c}
B(X_1, X_2, X_3, X_4) + B(X_2, X_1, X_3, X_4) = 0,\\
B(X_1, X_2, X_3, X_4) - B(X_3, X_4, X_1, X_2) = 0,\\
B(X_1, X_2, X_3, X_4) + B(X_2, X_3, X_1, X_4) + B(X_3, X_1, X_2, X_4) = 0.
\end{array}\right.
\eea
Solving the above equations we get the result.\\
\indent Thus if $B$ is a generalized curvature tensor then $B$ can be written as
\be\label{eqgen}
B = b_0 R + b_1 g\wedge S + b_2 r g\wedge g,
\ee
where $b_0$, $b_1$ and $b_2$ are scalars.\\
\indent We note that the equation $B(X_1, X_2, X_3, X_4) + B(X_2, X_3, X_1, X_4) + B(X_3, X_1, X_2, X_4) = 0$ can be omitted from the system of equations (\ref{eq3.1}) keeping the solution unaltered. Thus the tensor $B$ turns out to be a generalized curvature tensor if and only if
\beb
&B(X_1, X_2, X_3, X_4) + B(X_2, X_1, X_3, X_4) = 0,&\\
&B(X_1, X_2, X_3, X_4) - B(X_3, X_4, X_1, X_2) = 0.&
\eeb
\begin{lem}\label{lem3.12}
The tensor $B$ is a proper generalized curvature tensor if and only if $B$ is some constant multiple of $R$.
\end{lem}
\noindent\textbf{Proof:} Let $B$ be a proper generalized curvature tensor. Then obviously $B$ is a generalized curvature tensor and hence it can be written as
$$B = b_0 R + b_1 g\wedge S + b_2 r g\wedge g,$$
where $b_0$, $b_1$ and $b_2$ are scalars. Now $(g\wedge S)$ and $r (g\wedge g)$ both are not proper generalized curvature tensors. Hence for the tensor $B$ to be proper generalized curvature tensor, the scalars $b_1$ and $b_2$ must be zero (since $R$, $(g\wedge S)$ and $r (g\wedge g)$ are independent). Then $B=a_0 R$. Now from condition of proper generalized curvature tensor , we get $b_0 =$ constant. This completes the proof.
\begin{lem}\label{lem3.13}
The endomorphism operator $\mathscr B(X,Y)$ is skew-symmetric if $$a_2 = - a_6, \ a_3 = - a_5, \  a_8 = - a_9, \ a_1 = a_4 = a_7 = a_{10} = 0.$$
\end{lem}
\noindent\textbf{Proof:} From the Lemma \ref{lem5.22}, the operator $\mathscr B$ is skew-symmetric if
$$B(X_1,X_2,X_3,X_4) = -B(X_1,X_2,X_4,X_3) \mbox{ for all } X_1,X_2,X_3,X_4.$$
Thus the result follows from the solution of the equation $B(X_1, X_2, X_3, X_4) + B(X_1, X_2, X_4, X_3) = 0$.
%
\section{\bf{Main Results}}\label{Main}
In this section we first classify the set $\mathscr B$ of all $B$-tensors with respect to the contraction and then find out the equivalency 
of some structures for each class members. This classification can express in tree diagram as follows:\\

\vspace{0.1in}
$\put(115,0){\framebox(130,30){Class of $B$-tensors ($\mathscr B$)}}
\put(180,-1){\vector(0,-1){15}}
\put(80,-16){\line(1,0){200}}
\put(80,-16){\vector(0,-1){15}}    \put(280,-16){\vector(0,-1){15}}
\put(30,-61){\framebox(100,30){All $(^{ij}S) = 0$}}
\put(230,-61){\framebox(100,30){Some $(^{ij}S) \neq 0$}}
\put(280,-62){\vector(0,-1){15}}
\put(100,-77){\line(1,0){240}}
\put(100,-77){\vector(0,-1){15}}    \put(340,-77){\vector(0,-1){15}}
\put(50,-122){\framebox(100,30){All $(^{ij}p) = 0$}}
\put(290,-122){\framebox(100,30){Some $(^{ij}p) \neq 0$}}
\put(340,-122){\vector(0,-1){15}}
\put(160,-137){\line(1,0){240}}
\put(160,-137){\vector(0,-1){15}}    \put(400,-137){\vector(0,-1){15}}
\put(110,-182){\framebox(100,30){All $(^{ij}r) = 0$}}
\put(350,-182){\framebox(100,30){Some $(^{ij}r) \neq 0$}}\\
$
\vspace{0.1in}

Thus we get four different classes of $B$-tensors with respect to contraction given as follows:\\
(i) {\bf \textit{Class 1}:} In this class $(^{ij}S) = 0$ for all $i,j \in \{1,2,3,4\}$. Then we get dependency of $a_i$'s for this class as
\bea\label{eq4.1}
&&\left\{\begin{array}{c}
a_0 = -a_9(n-2)(n-1),\  a_1 = a_7(n-2),\\
a_2 = a_5 = -a_7 + (n-1)a_9,\  a_3 = a_6 = -(n-1)a_9,\  a_4 = a_7,\\
a_8 = -a_9 + \frac{a_7}{(n-1)},\   a_{10} = -\frac{a_7}{(n-1)}.
\end{array}\right.
\eea
\indent An example of such class of $B$-tensors is conformal curvature tensor $C$. We take $C$ as the representative member of this class.\\
(ii) {\bf \textit{Class 2}:} In this class $(^{ij}S) \neq 0$ for some $i,j \in \{1,2,3,4\}$ but $(^{ij}p) = 0$ for all $i,j \in \{1,2,3,4\}$. We get the dependency of $a_i$'s for this class that (\ref{eq4.1}) does not satisfy (i.e. one of $a_4 + a_8 + a_9 + n a_{10},\ a_3 + a_8 + n a_9 + a_{10},\ a_5 + n a_8 + a_9 + a_{10},\ a_2 + n a_8 + a_9 + a_{10},\ a_6 + a_8 + n a_9 + a_{10},\ a_7 + a_8 + a_9 + n a_{10}$ is non-zero) but
\bea\label{eq4.2}
&&a_0 = a_6(n-2),\ a_1 = a_7(n-2),\ a_2 = a_5 = -a_6 - a_7,\ a_3 = a_6,\ a_4 = a_7.
\eea
\indent An example of such class of $B$-tensors is conharmonic curvature tensor $K$. We take $K$ as the representative member of this class.\\
(iii) {\bf \textit{Class 3}:} In this class $(^{ij}p) \neq 0$ for some $i,j \in \{1,2,3,4\}$ but $(^{ij}r) = 0$ for all $i,j \in \{1,2,3,4\}$. Then for this class $a_i$'s does not satisfy (\ref{eq4.1}) but
\bea\label{eq4.3}
&&\left\{\begin{array}{c}
a_0 = (n-1)(a_3 + a_6 + n a_9), \ a_8 = -\frac{a_1 + (n-1)(a_2 + a_3 + a_5 + a_6 + n a_9)}{n(n-1)},\\
a_{10} = \frac{n(a_1 - (n-1)(a_4 + a_7))}{n(n-1)}.
\end{array}\right.
\eea
\indent Examples of such class of $B$-tensors are $W$, $P$, $\mathcal M$, $P^*$, $\mathcal W_0$, $\mathcal W_1$, $\mathcal W^*_3$. We take $W$ as the representative member of this class. We note that in this case $(^{ij}p) + n (^{ij}q) = 0$.\\
(iv) {\bf \textit{Class 4}:} In this class $(^{ij}p) \mbox{ and } (^{kl}r) \neq 0$ for some $i,j,k,l \in \{1,2,3,4\}$. For this class $a_i$'s does not satisfy (\ref{eq4.2}) and (\ref{eq4.3}).\\
\indent Examples of such class of $B$-tensors are $R$, $\mathcal W_0^*$, $\mathcal W_1^*$, $\mathcal W_2$, $\mathcal W_2^*$, $\mathcal W_3$, $\mathcal W_4$, $\mathcal W_4^*$, $\mathcal W_5$, $\mathcal W_5^*$, $\mathcal W_6$, $\mathcal W_6^*$, $\mathcal W_7$, $\mathcal W_7^*$, $\mathcal W_8$, $\mathcal W_8^*$, $\mathcal W_9$, $\mathcal W_9^*$. We take $R$ as the representative member of this class.\\
\indent We first discuss about the linear combination of $B$-tensors over $C^{\infty}(M)$. We consider two $B$-tensors $\bar B$ and $\tilde B$ with their $(i$-$j)$-th contraction tensors $(^{ij}\bar p) S + (^{ij}\bar q) r g$ and $(^{ij}\tilde p) S + (^{ij}\tilde q) r g$ respectively. Now consider a linear combination $\acute{B} = \mu \bar B + \eta \tilde B$ of $\bar B$ and $\tilde B$, where $\mu$ and $\eta$ are two scalars. Then 
$\acute B$ is a $B$-tensor with $(i$-$j)$-th contraction tensor $(^{ij}\acute p) S + (^{ij}\acute q) r g$, where $(^{ij}\acute p) = (^{ij}\bar p) + (^{ij}\tilde p)$ and $(^{ij}\acute q) = (^{ij}\bar q) + (^{ij}\tilde q)$.
 Now if both $\bar B$ and $\tilde B$ belong to class 1, then $(^{ij}\bar p) = (^{ij}\bar q) = (^{ij}\tilde p) = (^{ij}\tilde q) =0$ and thus $(^{ij}\acute p) = (^{ij}\acute q) =0$, i.e., $\acute B$ also remains a member of class 1. So class 1 is closed under linear combination over $C^{\infty}(M)$.
 If both $\bar B$ and $\tilde B$ belong to class 2, then $(^{ij}\bar p) = (^{ij}\tilde p) = 0$ for all $i,j$ but $(^{ij}\bar q)$ and $(^{ij}\tilde q)$ are not zero for all $i,j$ and thus $(^{ij}\acute p) = 0$. Now if $(^{ij}\acute q) = 0$ for all $i,j$, then $\acute B$ belongs to class 1, otherwise it remains a member of class 2. Here the condition $(^{ij}\acute q) = 0$ can be expressed explicitly as
\bea\label{c1}
&\left(\bar a_8+\bar a_9+\bar a_{10}\right) \mu+\left(\tilde a_8+\tilde a_9+\tilde a_{10}\right) \eta =0,&\\\nonumber
&(\mu\bar a_0 +\eta\tilde a_0) + (n-1)(n-2)(\mu\bar a_9 +\eta\tilde a_9) = 0,&\\\nonumber
&(\mu\bar a_2 +\eta\tilde a_2) = (n-1)\left[\mu(\bar a_9+\bar a_{10}) +\eta(\tilde a_9+\tilde a_{10})\right].&
\eea
Again if both $\bar B$ and $\tilde B$ belong to class 3, then $(^{ij}\bar p) + n(^{ij}\bar q)  = (^{ij}\tilde p) + n(^{ij}\tilde q) = 0$ for all $i,j$ but $(^{ij}\bar p)$ and $(^{ij}\tilde p)$ are not zero for all $i,j$ and thus $(^{ij}\acute p) + n(^{ij}\acute q) = 0$. Now if $(^{ij}\acute p) = 0$ or $(^{ij}\acute q) = 0$ for all $i,j$, then $\acute B$ belongs to class 1, otherwise it remains a member of class 3. Here the condition $(^{ij}\acute p) = 0$ and $(^{ij}\acute q) = 0$ are same and can be expressed explicitly as
\bea\label{c2}
&\left(\bar a_8+\bar a_9+\bar a_{10}\right) \mu+\left(\tilde a_8+\tilde a_9+\tilde a_{10}\right) \eta =0,&\\\nonumber
&(\mu\bar a_0 +\eta\tilde a_0) = (\mu\bar a_2 +\eta\tilde a_2) = (\mu\bar a_3 +\eta\tilde a_3) = (\mu\bar a_5 +\eta\tilde a_5) = -(n-1)(n-2)(\mu\bar a_9 +\eta\tilde a_9),&\\\nonumber
&(\mu\bar a_4 +\eta\tilde a_4) = (n-1)\left[\mu(\bar a_8+\bar a_9) +\eta(\tilde a_8+\tilde a_9)\right].&
\eea
We also note that if both $\bar B$ and $\tilde B$ belong to class 4, then $\acute B$ belongs to any one of the class according as their defining condition. Now if $\bar B$ belongs to class 1 then $\acute B$ belongs to the class as that of $\tilde B$. If $\bar B$ belongs to class 4 then $\acute B$ is of class 1 whether $\tilde B$ may belongs to class 2 or 3. Again if $\bar B$ belongs to class 2 and $\tilde B$ belongs to class 3, then obviously $\acute B$ becomes a member of class 4. Thus we can state the following:
\begin{thm}
$(i)$ Linear combinations of any two members of class 1 over $C^{\infty}(M)$ are the members of class 1.\\
$(ii)$ Linear combinations of any two members of class 2 over $C^{\infty}(M)$ are the members of class 1 or class 2 according as (\ref{c1}) holds or does not hold.\\
$(iii)$ Linear combinations of any two members of class 3 over $C^{\infty}(M)$ are the members of class 1 or class 3 according as (\ref{c2}) holds or does not hold.\\
$(iv)$ Linear combinations of any two members of class 4 over $C^{\infty}(M)$ may belong to any class among the four classes.\\
$(v)$ Linear combinations of any member of class 1 with any other member of any one of the remaining three classes over $C^{\infty}(M)$ belongs to the latter class.\\
$(vi)$ Linear combinations of any member of class 4 with any other member of any one of the remaining three classes over $C^{\infty}(M)$ belongs to class 4.\\
$(vi)$ Linear combinations of any member of class 2 with any other member of class 3 over $C^{\infty}(M)$ is a member of class 4.
\end{thm}
\indent From above we note that the tensor $B$ belongs to any one of the four classes according to the dependency of the coefficients $a_i$'s. The $\mathcal T$-curvature tensor and $C'$ may also belongs to any class. The curvature tensor $C^*$, $\widetilde W$ and $W^*$ are combination of two or more other $B$-tensors and thus they may belong to more than one class according as the coefficients of such combinations. Thus $C^*$ is a member of class 3 if $a_0 + (n-2)a_2 \neq 0$, otherwise it reduces to the conformal curvature tensor and becomes a member of class 1. Again, $\widetilde W$ is a member of class 3 if $a_0 -a_2 + (n-1)a_5 \neq 0$, otherwise it belongs to class 1. And $W^*$ is a member of class 4 if $b \neq 0$, otherwise it belongs to class 3. We also note that $P^*$ is the combination of $P$ and $W$, both of them are in class 3 and $P^*$ remains a member of class 3.\\
\indent We now discuss the equivalency of flatness, symmetry type, recurrent type, semisymmetry type and other various curvature restrictions for the above four classes of $B$-tensors.
\begin{thm}\label{thm4.1}
Flatness of all $B$-tensors of any class among the classes $1$- $4$ are equivalent to the flatness of the representative member of that class. (Flatness of all $B$-tensors of each class are equivalent.)
\end{thm}
%
\noindent\textbf{Proof:} We first consider that the tensor $B$ belongs to class 1 i.e. $a_i$'s satisfy (\ref{eq4.1}). We have to show $B = 0$ if and only if $C = 0$. Now
\bea\label{e1}
\nonumber B_{ijkl}-(a_0 C_{ijkl} + a_1 C_{ikjl}) &=& \frac{a_0 \left(-g_{jl}S_{ik}+g_{jk} S_{il}+g_{il} S_{jk}-g_{ik} S_{jl}\right)}{(n-2)}-\frac{a_0 \left(g_{il} g_{jk}-g_{ik} g_{jl}\right) r}{(n-2) (n-1)}\\
&&+\frac{a_1 \left(-g_{kl} S_{ij}+g_{jk} S_{il}+g_{il} S_{jk}-g_{ij} S_{kl}\right)}{n-2}-\frac{a_1 \left(g_{il} g_{jk}-g_{ij} g_{kl}\right) r}{(n-1) (n-2)}\\
\nonumber &&+a_2 g_{il} S_{jk}+a_3 g_{jl} S_{ik}+a_4 g_{kl} S_{ij}+a_5 g_{jk} S_{il}+a_6 g_{ik} S_{jl}+a_7 g_{ij} S_{kl}\\
\nonumber &&+r \left(a_8 g_{il} g_{jk}+a_9 g_{ik} g_{jl}+a_{10} g_{ij} g_{kl}\right).
\eea
As $B$ is of class 1 so simplifying the above and using (\ref{eq4.1}) we get $B_{ijkl} - (a_0 C_{ijkl} + a_1 C_{ikjl}) = 0$. Again $(a_0 C_{ijkl} + a_1 C_{ikjl}) = 0$ if and only if $C_{ijij} = 0$, i.e. if and only if $C_{ijkl} = 0$ (by Lemma \ref{lem3.6}). Thus we get $B=0$ if and only if $C = 0$.\\
\indent Next we consider that the tensor $B$ belongs to class 2 i.e. $a_i$'s satisfy (\ref{eq4.2}) but not (\ref{eq4.1}). We have to show $B = 0$ if and only if $K = 0$. Now
\beb
B_{ijkl} -(a_0 K_{ijkl} + a_1 K_{ikjl}) &=& \frac{a_0 \left(-g_{jl}S_{ik}+g_{jk} S_{il}+g_{il} S_{jk}-g_{ik} S_{jl}\right)}{n-2}\\
&+&\frac{a_1 \left(-g_{kl} S_{ij}+g_{jk} S_{il}+g_{il} S_{jk}-g_{ij} S_{kl}\right)}{n-2}\\
&+&a_2 g_{il} S_{jk}+a_3 g_{jl} S_{ik}+a_4 g_{kl} S_{ij}+a_5 g_{jk} S_{il}+a_6 g_{ik} S_{jl}+a_7 g_{ij} S_{kl}\\
&+&r \left(a_8 g_{il} g_{jk}+a_9 g_{ik} g_{jl}+a_{10} g_{ij} g_{kl}\right).
\eeb
As $B$ is of class 2 so simplifying the above and using (\ref{eq4.2}) we get
\be\label{e2}
B_{ijkl} -(a_0 K_{ijkl} + a_1 K_{ikjl}) = r \left(a_8 g_{il} g_{jk}+a_9 g_{ik} g_{jl}+a_{10} g_{ij} g_{kl}\right).
\ee
Now as $B$ and $K$ are both of class 2 so vanishing of any one of $B$ or $K$ implies $r = 0$ and then
$$B_{ijkl} -(a_0 K_{ijkl} + a_1 K_{ikjl}) = 0.$$
Again, $(a_0 K_{ijkl} + a_1 K_{ikjl}) = 0$ if and only if $K_{ijij} = 0$, i.e. if and only if $K_{ijkl} = 0$ (by Lemma \ref{lem3.6}). Thus from (\ref{e2}), we get $B=0$ if and only if $K = 0$.\\
\indent Again we consider that the tensor $B$ belongs to class 3 i.e. $a_i$'s satisfies (\ref{eq4.3}) but not (\ref{cor4.1}). We have to show $B = 0$ if and only if $W = 0$. Now
\beb
B_{ijkl} -(a_0 W_{ijkl} + a_1 W_{ikjl}) &=&
\frac{a_0 \left(g_{jk} g_{il}-g_{jl}g_{ik}+g_{il} g_{jk}-g_{ik} g_{jl}\right)}{n(n-i)}\\
&+&\frac{a_1 \left(g_{jk} g_{il}-g_{kl} g_{ij}+g_{il} g_{jk}-g_{ij} g_{kl}\right)}{n(n-i)}\\
&+&a_2 g_{il} S_{jk}+a_3 g_{jl} S_{ik}+a_4 g_{kl} S_{ij}+a_5 g_{jk} S_{il}+a_6 g_{ik} S_{jl}+a_7 g_{ij} S_{kl}\\
&+&r \left(a_8 g_{il} g_{jk}+a_9 g_{ik} g_{jl}+a_{10} g_{ij} g_{kl}\right).
\eeb
As $B$ is of class 3 so simplifying the above and using (\ref{eq4.3}) we get
\bea\label{e3}
\nonumber B_{ijkl} -(a_0 W_{ijkl} + a_1 W_{ikjl}) &=&
\frac{r}{n} [\left(a_2 g_{il} g_{jk}+a_3 g_{ik} g_{jl}+a_4 g_{ij} g_{kl}+a_5 g_{il} g_{jk}+a_6 g_{ik} g_{jl}+a_7 g_{ij} g_{kl}\right)]\\
&-& [a_2 g_{il} S_{jk}+a_3 g_{jl} S_{ik}+a_4 g_{kl} S_{ij}+a_5 g_{jk} S_{il}+a_6 g_{ik} S_{jl}+a_7 g_{ij} S_{kl}].
\eea
Now as $B$ and $W$ are both of class 3 so vanishing of any one of $B$ or $W$ implies $S = \frac{r}{n} g$ and then
$$\nonumber B_{ijkl} -(a_0 W_{ijkl} + a_1 W_{ikjl}) = 0.$$
Again, $(a_0 W_{ijkl} + a_1 W_{ikjl}) = 0$ if and only if $W_{ijij} = 0$, i.e., if and only if $W_{ijkl} = 0$ (by Lemma \ref{lem3.6}). Thus from (\ref{e3}), we get $B=0$ if and only if $W = 0$.\\
\indent Finally, we consider that the tensor $B$ belongs to class 4. We have to show $B = 0$ if and only if $R = 0$. Now as $B$ and $R$ are both of class 4 so vanishing of any one of $B$ or $R$ implies $S = 0$. Thus by the Lemma \ref{lem3.8} we can conclude that $B = 0$ if and only if $R = 0$.
This completes the proof.\\
\indent From the proof of the above we can state the following:
\begin{cor}\label{cor4.1}
If $B$ belongs to any class out of the above four classes then $R = 0$ implies $B = 0$ and $B = 0$ implies $C = 0$.
\end{cor}
\noindent \textbf{Proof:} We know that $R=0$ implies $S= 0$ and $r=0$. So the first part of the proof is obvious. Again we know that $R=0$ or $W=0$ or $K=0$ all individually implies $C=0$. So by above theorem the proof of the second part is done.\\
\indent We now discuss the above four classes of $B$-tensors as equivalence classes of an equivalence relation on $\mathscr B$, the set of all $B$-tensors. Consider a relation $\rho$ on $\mathscr B$ defined by $B_1 \rho B_2$ if and only if $B_1$-flat (i.e. $B_1 = 0$) $\Leftrightarrow$ $B_2$-flat (i.e. $B_2 = 0$) , for all $B_1, B_2 \in \mathscr B$. It can be easily shown that $\rho$ is an equivalence relation. We conclude from Theorem \ref{thm4.1} that all $B$-tensors of class 1 are related to the conformal curvature tensor $C$, all $B$-tensors of class 2 are related to the conharmonic curvature tensor $K$, all $B$-tensors of class 3 are related to the concircular curvature tensor $W$, all $B$-tensors of class 4 are related to the Riemann-Christoffel curvature tensor $R$. Thus class 1 is the  $\rho$-equivalence class [$C$], class 2 is the $\rho$-equivalence class [$K$], class 3 is the $\rho$-equivalence class [$W$] and class 4 is the $\rho$-equivalence class [$R$].
\begin{thm}\label{cl-1}
$($Characteristic of class $1)$ $(i)$ All tensors of class $1$ are of the form
				$$a_0 C_{ijkl} + a_1 C_{ikjl}$$
and the only generalized curvature tensor of this class is conformal curvature tensor upto a scalar multiple.\\
$(ii)$ All curvature restrictions of type $1$ on any $B$-tensor of class $1$ are equivalent, if $a_0$ and $a_1$ are constants.\\
$(iii)$ All curvature restrictions of type $2$ on any $B$-tensor of class $1$ are equivalent.
\end{thm}
\noindent \textbf{Proof:} We see that if $B$ is of class 1, then from (\ref{e1}), $B_{ijkl} = \left[a_0 C_{ijkl} + a_1 C_{ikjl}\right]$. Now for $B$ to be a generalized curvature tensor (\ref{eq5.3}) fulfilled and we get the form of generalized curvature tensor of this class.\\
Again consider any curvature restriction by an operator $\mathcal L$ on $B$, we have $\mathcal L B = 0$, which implies $\mathcal L\left[a_0 C_{ijkl} + a_1 C_{ikjl}\right] = 0$. Then by Lemma \ref{lem3.6} we get the result.
\begin{thm}\label{cl-2}
$($Characteristic of class $2)$ $(i)$ All tensors of class $2$ are of the form
				$$a_0 K_{ijkl} + a_1 K_{ikjl} + r \left(a_8 g_{il} g_{jk}+a_9 g_{ik} g_{jl}+a_{10} g_{ij} g_{kl}\right)$$
such that $a_8 = \frac{a0+a1}{(n-2) (n-1)}, a_9 = -\frac{a0}{(n-2) (n-1)}, a_{10} = -\frac{a1}{(n-2) (n-1)}$ are not satisfy all together, otherwise it belongs to class $1$. The generalized curvature tensor of this class are of the form $a_0 K + a_8 r G$, $a_8 \neq \frac{a0}{(n-1) (n-2)}$.\\
$(ii)$ All commutative curvature restrictions of type $1$ on any $B$-tensor of class $2$ are equivalent, if $a_0$, $a_1$, $a_8$, $a_9$ and $a_{10}$ are constants.\\
$(iii)$ All commutative curvature restrictions of type $2$ on any $B$-tensor of class $2$ are equivalent.
\end{thm}
\noindent \textbf{Proof:} We see that if $B$ is of class 2, then from (\ref{e2}),
$$B_{ijkl} = a_0 K_{ijkl} + a_1 K_{ikjl} + r \left(a_8 g_{il} g_{jk}+a_9 g_{ik} g_{jl}+a_{10} g_{ij} g_{kl}\right).$$
Now for $B$ to be generalized curvature tensor, (\ref{eq5.3}) fulfilled and we get the form of generalized curvature tensor of this class as required.\\
Again consider any curvature restriction by a commutative operator $\mathcal L$ on $B$ i,e., $\mathcal L B = 0$, which implies
$$\mathcal L \left[a_0 K_{ijkl} + a_1 K_{ikjl} + r \left(a_8 g_{il} g_{jk}+a_9 g_{ik} g_{jl}+a_{10} g_{ij} g_{kl}\right)\right] = 0.$$
Now if $\mathcal L$ is of first type and $a_0$, $a_1$, $a_8$, $a_9$ and $a_{10}$ are constants, then $\mathcal L[a_0 K_{ijkl} + a_1 K_{ikjl}]=0$ and if $\mathcal L$ is of second type then automatically $\mathcal L[a_0 K_{ijkl} + a_1 K_{ikjl}]=0$. Thus by Lemma \ref{lem3.6} we get the result.
\begin{thm}\label{cl-3}
$($Characteristic of class $3)$ $(i)$ All tensors of class $3$ are of the form
\beb
&&(a_0 W_{ijkl} + a_1 W_{ikjl}) +\left[a_2 g_{il} S_{jk} + a_3 g_{jl} S_{ik} + a_4 g_{kl} S_{ij} + a_5 g_{jk} S_{il} + a_6 g_{ik} S_{jl} + a_7 g_{ij} S_{kl}\right]\\
&& - \frac{r}{n}\left[(a_2+a_5) g_{il} g_{jk} + (a_3+a_6) g_{jl} g_{ik} + (a_4+a_7) g_{kl} g_{ij}\right]
\eeb
such that $a_2 = a_5 = -\frac{a0 + a1}{n-2}, a_3 = a_6 = \frac{a0}{n-2}, a_4 = a_7 = \frac{a1}{n-2}$ are not satisfy all together, otherwise it belongs to class $1$. The generalized curvature tensor of this class are of the form $a_0 W + a_2 \left[g\wedge S - \frac{r}{n}G\right]$, $a_2 \neq \frac{a0}{(n-1)(n-2)}$.\\
$(ii)$ All commutative curvature restrictions of type $1$ on any $B$-tensor of class $3$ are equivalent, if $a_0$, $a_1$, $a_2$, $a_3$, $a_4$, $a_5$, $a_6$ and $a_7$ are constants.\\
$(iii)$ All commutative curvature restrictions of type $2$ on any $B$-tensor of class $3$ are equivalent.
\end{thm}
\noindent \textbf{Proof:} The proof is similar to the proof of the Theorem \ref{cl-2}.
\begin{thm}\label{cl-4}
$($Characteristic of class $4)$ $(i)$ All commutative curvature restrictions of type $1$ on any $B$-tensor of class $4$ are equivalent, if $a_i$'s are all constants.\\
$(ii)$ All commutative curvature restrictions of type $2$ on any $B$-tensor of class $4$ are equivalent.
\end{thm}
\noindent \textbf{Proof:} Consider a commutative operator $\mathcal L$ and $B$ is of class 1, such that $\mathcal L B = 0$. Now if $\mathcal L$ is commutative 1st type and all $a_i$'s are constant, then by taking contraction we get $\mathcal L S = 0$ and $\mathcal L (r) = 0$ as $a_i$'s are all constant. Putting these in the expression of $\mathcal L B = 0$ we get, $\mathcal L R = 0$. Again if  $\mathcal L$ is of commutative 2nd type, then contraction yields $\mathcal L S = 0$ and $\mathcal L(r) = 0$. Substituting these in the expression of $\mathcal L B = 0$, we get $\mathcal L R = 0$. This complete the proof.\\
\indent From the proofs of the above four characteristic theorem as similar of Corollary \ref{cor4.1}, we can state the following:
\begin{cor}
If the tensor $B$ belongs to any one of the class and $\mathcal L$ is a commutative operator. If $\mathcal L$ is of type 2, then\\
(i) $\mathcal L R = 0$ implies $\mathcal L B = 0$ and\\
(ii) $\mathcal L B = 0$ implies $\mathcal L C = 0$.\\
The results are also true for the case of type 1 if the coefficients of $B$ are all constant.
\end{cor}
\noindent \textbf{Proof:} First consider the case $\mathcal L$ to be of 2nd type and commutative. Then $\mathcal L R=0$ implies $\mathcal L S= 0$ and $\mathcal L r=0$. So the first part of the proof is obvious. Again we can easily check that $\mathcal L R=0$ or $\mathcal L W=0$ or $\mathcal L K=0$ all individually implies $\mathcal L C=0$. So by the above four characteristic theorem the proof of the second part is done.\\ The proof for the next case (i.e., $\mathcal L$ to be of type 1 and commutative with the coefficients of $B$'s are all constant) is similar to above.\\
%
\indent We now state some results which will be used to show the coincidence of class 3 ad class 4 for the symmetry and recurrency condition.
\begin{lem}\label{lem5.1}\cite{RT} 
Locally symmetric and projectively symmetric semi-Riemannian manifolds are equivalent.
\end{lem}
\begin{lem}\label{lem5.2}$($\cite{Glodek}, \cite{ol}, \cite{mik}, \cite{mik2}$)$ 
Every concircularly recurrent as well as projective recurrent manifold is necessarily a recurrent manifold with the same recurrence form.
\end{lem}
From the above four Characteristic theorems of the four classes and the Lemma \ref{lem5.1} and \ref{lem5.2} we can state the results expressed in a table for 1st type operator such that in the following table all condition(s) in a block are equivalent.
\begin{table}[H]
\begin{center}
{\footnotesize\begin{tabular}{|c|c|c|c|}\hline
\multicolumn{1}{|c|}{\textbf{Class 1}}  &  \multicolumn{1}{|c|}{\textbf{Class 2}}  &  \multicolumn{1}{|c|}{\textbf{Class3}}  &  \multicolumn{1}{|c|}{\textbf{\textbf{Class4}}}\\\hline

\multicolumn{1}{|c|}{}&
\multicolumn{1}{|c|}{}&
\multicolumn{1}{|c|}{$W = 0$, $P = 0$, $\mathcal M = 0$,}&
\multicolumn{1}{|c|}{$R = 0$, $\mathcal W_0^* = 0$, $\mathcal W_1^* = 0$, $\mathcal W_2 = 0$,}\\
\multicolumn{1}{|c|}{$C = 0$}&
\multicolumn{1}{|c|}{$K = 0$}&
\multicolumn{1}{|c|}{$P^* = 0$, $\mathcal W_0 = 0$,}&
\multicolumn{1}{|c|}{$\mathcal W_2^* = 0$, $\mathcal W_3 = 0$,}\\
\multicolumn{1}{|c|}{       }&
\multicolumn{1}{|c|}{       }&
\multicolumn{1}{|c|}{$\mathcal W_1 = 0$, $\mathcal W_3^* = 0$}&
\multicolumn{1}{|c|}{$\mathcal W_i = 0$, $\mathcal W_i^* = 0$, for all $i = 4,5,\cdots 9$}\\\hline

\multicolumn{1}{|c|}{}&
\multicolumn{1}{|c|}{}&
\multicolumn{1}{|r}{$\nabla W = 0$, $\nabla P = 0$, $\nabla \mathcal M = 0$,}&
\multicolumn{1}{l|}{$\nabla R = 0$, $\nabla \mathcal W_0^* = 0$, $\nabla \mathcal W_1^* = 0$, $\nabla \mathcal W_2 = 0$,}\\
\multicolumn{1}{|c|}{$\nabla C = 0$}&
\multicolumn{1}{|c|}{$\nabla K = 0$}&
\multicolumn{1}{|r}{$\nabla P^* = 0$, $\nabla \mathcal W_0 = 0$,}&
\multicolumn{1}{l|}{$\nabla \mathcal W_2^* = 0$, $\nabla \mathcal W_3 = 0$,}\\
\multicolumn{1}{|c|}{       }&
\multicolumn{1}{|c|}{       }&
\multicolumn{1}{|r}{$\nabla \mathcal W_1 = 0$, $\nabla \mathcal W_3^* = 0$,}&
\multicolumn{1}{l|}{$\nabla \mathcal W_i = 0$, $\nabla \mathcal W_i^* = 0$, for all $i = 4,5,\cdots 9$}\\\hline

\multicolumn{1}{|c|}{}&
\multicolumn{1}{|c|}{}&
\multicolumn{1}{|r}{$\kappa W = 0$, $\kappa P = 0$, $\kappa \mathcal M = 0$,}&
\multicolumn{1}{l|}{$\kappa R = 0$, $\kappa \mathcal W_0^* = 0$, $\kappa \mathcal W_1^* = 0$, $\kappa \mathcal W_2 = 0$,}\\
\multicolumn{1}{|c|}{$\kappa C = 0$}&
\multicolumn{1}{|c|}{$\kappa K = 0$}&
\multicolumn{1}{|r}{$\kappa P^* = 0$, $\kappa \mathcal W_0 = 0$,}&
\multicolumn{1}{l|}{$\kappa \mathcal W_2^* = 0$, $\kappa \mathcal W_3 = 0$,}\\
\multicolumn{1}{|c|}{       }&
\multicolumn{1}{|c|}{       }&
\multicolumn{1}{|r}{$\kappa \mathcal W_1 = 0$, $\kappa \mathcal W_3^* = 0$,}&
\multicolumn{1}{l|}{$\kappa \mathcal W_i = 0$, $\kappa \mathcal W_i^* = 0$, for all $i = 4,5,\cdots 9$}\\\hline

\multicolumn{1}{|c|}{}&
\multicolumn{1}{|c|}{}&
\multicolumn{1}{|c|}{$L^s W = 0$, $L^s P = 0$,}&
\multicolumn{1}{|c|}{$L^s R = 0$, $L^s \mathcal W_0^* = 0$, $L^s \mathcal W_1^* = 0$,}\\
\multicolumn{1}{|c|}{$L^s C = 0$}&
\multicolumn{1}{|c|}{$L^s K = 0$}&
\multicolumn{1}{|c|}{$L^s \mathcal M = 0$, $L^s P^* = 0$,}&
\multicolumn{1}{|c|}{$L^s \mathcal W_2 = 0$, $L^s \mathcal W_2^* = 0$, $L^s \mathcal W_3 = 0$,}\\
\multicolumn{1}{|c|}{       }&
\multicolumn{1}{|c|}{       }&
\multicolumn{1}{|c|}{$L^s \mathcal W_0 = 0$, $L^s \mathcal W_1 = 0$,}&
\multicolumn{1}{|c|}{$L^s \mathcal W_i = 0$, $L^s \mathcal W_i^* = 0$,}\\
\multicolumn{1}{|c|}{       }&
\multicolumn{1}{|c|}{       }&
\multicolumn{1}{|c|}{$L^s \mathcal W_3^* = 0$}&
\multicolumn{1}{|c|}{for all $i = 4,5,\cdots 9$}\\\hline

\multicolumn{1}{|c|}{}&
\multicolumn{1}{|c|}{}&
\multicolumn{1}{|c|}{$\kappa^s W = 0$, $\kappa^s P = 0$,}&
\multicolumn{1}{|c|}{$\kappa^s R = 0$, $\kappa^s \mathcal W_0^* = 0$, $\kappa^s \mathcal W_1^* = 0$,}\\
\multicolumn{1}{|c|}{$\kappa^s C = 0$}&
\multicolumn{1}{|c|}{$\kappa^s K = 0$}&
\multicolumn{1}{|c|}{$\kappa^s \mathcal M = 0$, $\kappa^s P^* = 0$,}&
\multicolumn{1}{|c|}{$\kappa^s \mathcal W_2 = 0$, $\kappa^s \mathcal W_2^* = 0$, $\kappa^s \mathcal W_3 = 0$,}\\
\multicolumn{1}{|c|}{       }&
\multicolumn{1}{|c|}{       }&
\multicolumn{1}{|c|}{$\kappa^s \mathcal W_0 = 0$, $\kappa^s \mathcal W_1 = 0$,}&
\multicolumn{1}{|c|}{$\kappa^s \mathcal W_i = 0$, $\kappa^s \mathcal W_i^* = 0$,}\\
\multicolumn{1}{|c|}{       }&
\multicolumn{1}{|c|}{       }&
\multicolumn{1}{|c|}{$\kappa^s \mathcal W_3^* = 0$}&
\multicolumn{1}{|c|}{for all $i = 4,5,\cdots 9$}\\\hline

\end{tabular}}
\end{center}
\caption{List of equivalent structures for 1st type operators}\label{stc-1}
\end{table}
\indent\indent We now prove a theorem for coincidence of class 1 with class 2 and coincidence of class 3 with class 4 for any commutative 2nd type operator.
\begin{thm}\label{thm5.6}
Let $\mathcal L$ be a commutative operator (i.e., $\mathcal L$ and contraction commute) of type $2$, Then the following holds:\\
(i) For any two $B$-tensor $B_1$ and $B_2$ of class $1$ and $2$ respectively, the conditions $\mathcal L B_1 = 0$ and $\mathcal L B_2 = 0$ are equivalent.\\
(ii) For any two $B$-tensor $B_1$ and $B_2$ of class $3$ and $4$ respectively, the conditions $\mathcal L B_1 = 0$ and $\mathcal L B_2 = 0$ are equivalent.
\end{thm}
\noindent \textbf{Proof:} From above four characteristic theorem it is clear that to prove this theorem it is sufficient to show that for a commutative 2nd type operator $\mathcal L$, $\mathcal L C = \mathcal L K$ and $\mathcal L W = \mathcal L R$. Now first consider $C$ and $K$. Then for any operator $\mathcal L$,
$$\mathcal L C = \mathcal L K + \mathcal L \left(\frac{r}{(n-1)(n-2)} G\right).$$
Thus if $\mathcal L$ is commutative 2nd type operator, then $\mathcal L \left(\frac{r}{(n-1)(n-2)} G\right) = 0$ and hence (i) is proved.\\
Again considering $R$ and $W$, we get
$$\mathcal L W = \mathcal L R + \mathcal L \left(\frac{r}{n(n-1)} G\right).$$
So if $\mathcal L$ is commutative 2nd type operator then $\mathcal L \left(\frac{r}{n(n-1)} G\right) = 0$ and hence (ii) is proved.\\
\indent From the above four Characteristic theorems of the classes and the Theorem \ref{thm5.6} we can state the results express in a table for 2nd type operator such that in the following table all conditions in a block of the table are equivalent.
\begin{table}[H]
\begin{center}
{\footnotesize
\begin{tabular}{|c|c|c|c|}\hline
\multicolumn{2}{|c|}{\textbf{Class 1 and Class 2}}&
\multicolumn{2}{|c|}{\textbf{Class 3 and Class 4}}\\\hline
\multicolumn{2}{|c|}{$R\cdot C=0$,}&
\multicolumn{2}{|c|}{$R\cdot W = 0$, $R\cdot P = 0$, $R\cdot P^* = 0$, $R\cdot \mathcal M = 0$, $R\cdot R = 0$,}\\
\multicolumn{2}{|c|}{$R\cdot K=0$}&
\multicolumn{2}{|c|}{$R\cdot \mathcal W_i = 0$, $R\cdot \mathcal W_i^* = 0$, for all $i = 0,1,\cdots 9$}\\\hline

\multicolumn{2}{|c|}{$C\cdot C=0$,}&
\multicolumn{2}{|c|}{$C\cdot W = 0$, $C\cdot P = 0$, $C\cdot P^* = 0$, $C\cdot \mathcal M = 0$, $C\cdot R = 0$,}\\
\multicolumn{2}{|c|}{$C\cdot K=0$}&
\multicolumn{2}{|c|}{$C\cdot \mathcal W_i = 0$, $C\cdot \mathcal W_i^* = 0$, for all $i = 0,1,\cdots 9$}\\\hline

\multicolumn{2}{|c|}{$K\cdot C=0$,}&
\multicolumn{2}{|c|}{$K\cdot W = 0$, $K\cdot P = 0$, $K\cdot P^* = 0$, $K\cdot \mathcal M = 0$, $K\cdot R = 0$,}\\
\multicolumn{2}{|c|}{$K\cdot K=0$}&
\multicolumn{2}{|c|}{$K\cdot \mathcal W_i = 0$, $K\cdot \mathcal W_i^* = 0$, for all $i = 0,1,\cdots 9$}\\\hline

\multicolumn{2}{|c|}{$W\cdot C=0$,}&
\multicolumn{2}{|c|}{$W\cdot W = 0$, $W\cdot P = 0$, $W\cdot P^* = 0$, $W\cdot \mathcal M = 0$, $W\cdot R = 0$,}\\
\multicolumn{2}{|c|}{$W\cdot K=0$}&
\multicolumn{2}{|c|}{$W\cdot \mathcal W_i = 0$, $W\cdot \mathcal W_i^* = 0$, for all $i = 0,1,\cdots 9$}\\\hline

\multicolumn{2}{|c|}{$R\cdot C=L Q(g,C)$,}&
\multicolumn{2}{|c|}{$R\cdot W = L Q(g,W)$, $R\cdot P = L Q(g,P)$, $R\cdot P^* = L Q(g,P^*)$, $R\cdot \mathcal M = L Q(g,\mathcal M)$,}\\
\multicolumn{2}{|c|}{$R\cdot K=L Q(g,K)$}&
\multicolumn{2}{|c|}{$R\cdot R = L Q(g,R)$, $R\cdot \mathcal W_i =L Q(g,\mathcal W_i)$, $R\cdot \mathcal W_i^* =L Q(g,\mathcal W_i^*)$, for all $i = 0,1,\cdots 9$}\\\hline

\multicolumn{2}{|c|}{$C\cdot C=L Q(g,C)$,}&
\multicolumn{2}{|c|}{$C\cdot W = L Q(g,W)$, $C\cdot P = L Q(g,P)$, $C\cdot P^* = L Q(g,P^*)$, $C\cdot \mathcal M = L Q(g,\mathcal M)$,}\\
\multicolumn{2}{|c|}{$C\cdot K=L Q(g,K)$}&
\multicolumn{2}{|c|}{$C\cdot R = L Q(g,R)$, $C\cdot \mathcal W_i =L Q(g,\mathcal W_i)$, $C\cdot \mathcal W_i^* =L Q(g,\mathcal W_i^*)$, for all $i = 0,1,\cdots 9$}\\\hline

\multicolumn{2}{|c|}{$W\cdot C=L Q(g,C)$,}&
\multicolumn{2}{|c|}{$W\cdot W = L Q(g,W)$, $W\cdot P = L Q(g,P)$, $W\cdot P^* = L Q(g,P^*)$, $W\cdot \mathcal M = L Q(g,\mathcal M)$,}\\
\multicolumn{2}{|c|}{$W\cdot K=L Q(g,K)$}&
\multicolumn{2}{|c|}{$W\cdot R = L Q(g,R)$, $W\cdot \mathcal W_i =L Q(g,\mathcal W_i)$, $W\cdot \mathcal W_i^* =L Q(g,\mathcal W_i^*)$, for all $i = 0,1,\cdots 9$}\\\hline

\multicolumn{2}{|c|}{$K\cdot C=L Q(g,C)$,}&
\multicolumn{2}{|c|}{$K\cdot W = L Q(g,W)$, $K\cdot P = L Q(g,P)$, $K\cdot P^* = L Q(g,P^*)$, $K\cdot \mathcal M = L Q(g,\mathcal M)$,}\\
\multicolumn{2}{|c|}{$K\cdot K=L Q(g,K)$}&
\multicolumn{2}{|c|}{$K\cdot R = L Q(g,R)$, $K\cdot \mathcal W_i =L Q(g,\mathcal W_i)$, $K\cdot \mathcal W_i^* =L Q(g,\mathcal W_i^*)$, for all $i = 0,1,\cdots 9$}\\\hline

\end{tabular}}
\end{center}
\caption{List of equivalent structures for 2nd type operators}\label{stc-2}
\end{table}
\indent Thus from the Table \ref{stc-2} we can state the following:
\begin{cor}\label{cor6.5}
$(1)$ The conditions $R\cdot R = 0$, $R\cdot W = 0$ and $R\cdot P = 0$ are equivalent.\\
$(2)$ The conditions $C\cdot R = 0$, $C\cdot W = 0$ and $C\cdot P = 0$ are equivalent.\\
$(3)$ The conditions $W\cdot R = 0$, $W\cdot W = 0$ and $W\cdot P = 0$ are equivalent.\\
$(4)$ The conditions $K\cdot R = 0$, $K\cdot W = 0$ and $K\cdot P = 0$ are equivalent.\\
$(5)$ The conditions $R\cdot C = 0$ and $R\cdot K = 0$ are equivalent.\\
$(6)$ The conditions $C\cdot C = 0$ and $C\cdot K = 0$ are equivalent.\\
$(7)$ The conditions $W\cdot C = 0$ and $W\cdot K = 0$ are equivalent.\\
$(8)$ The conditions $K\cdot C = 0$ and $K\cdot K = 0$ are equivalent.
\end{cor}
\indent It may be mentioned that the first four results of Corollary \ref{cor6.5} were proved in Theorem 3.3 of \cite{PV} in another way.
\begin{cor}
$(1)$ The conditions $R\cdot R = L_1 Q(g,R)$, $R\cdot W = L_1 Q(g,W)$ and $R\cdot P = L_1 Q(g,P)$ are equivalent.\\
$(2)$ The conditions $C\cdot R = L_2 Q(g,R)$, $C\cdot W = L_2 Q(g,W)$ and $C\cdot P = L_2 Q(g,P)$ are equivalent.\\
$(3)$ The conditions $W\cdot R = L_3 Q(g,R)$, $W\cdot W = L_3 Q(g,W)$ and $W\cdot P = L_3 Q(g,P)$ are equivalent.\\
$(4)$ The conditions $K\cdot R = L_4 Q(g,R)$, $K\cdot W = L_4 Q(g,W)$ and $K\cdot P = L_4 Q(g,P)$ are equivalent.\\
$(5)$ The conditions $R\cdot C = L_5 Q(g,C)$ and $R\cdot K = L_5 Q(g,K)$ are equivalent.\\
$(6)$ The conditions $C\cdot C = L_6 Q(g,C)$ and $C\cdot K = L_6 Q(g,K)$ are equivalent.\\
$(7)$ The conditions $W\cdot C = L_7 Q(g,C)$ and $W\cdot K = L_7 Q(g,K)$ are equivalent.\\
$(8)$ The conditions $K\cdot C = L_8 Q(g,C)$ and $K\cdot K = L_8 Q(g,K)$ are equivalent.\\
Here $L_i$, $(i=1,2,\cdots, 8)$ are scalars.
\end{cor}
%
We note that here the operator $\mathcal P$ for projective curvature tensor is not considered as $P$ is not skew-symmetric in 3rd and 4th places i.e. $\mathcal P$ is not commutative.
%
\section{\bf{Conclusion}}
%
Form the above discussion we see that the set $\mathscr B$ of all $B$-tensors can be partitioned into four equivalence classes [$C$] (or class 1), [$K$] (or class 2), [$W$] (or class 3) and [$R$] (or class 4) under the equivalence relation $\rho$ such that $B_1 \rho B_2$ holds if and only if $B_1 = 0$ implies $B_2 = 0$ and $B_2 = 0$ implies $B_1 = 0$, where $B_1, B_2 \in \mathscr B$. We conclude that\\
(i) study of any curvature restriction of type 1 (such as symmetric type, recurrent type, super generalized recurrent) on any $B$-tensor of class 1 with constant $a_i$'s is equivalent to the study of such type of curvature restriction on the conformal curvature tensor $C$ and also any curvature restriction of type 2 (such as semisymmetric type, pseudosymmetric type) on any $B$-tensor of class 2 is equivalent to the study of such type of curvature restriction on $C$. Thus for all such restrictions, each $B$-tensor of class 1 gives the same structure as that due to $C$.\\
(ii) study of a symmetric type and recurrent type curvature restrictions on any $B$-tensor of class 2 with constant $a_i$'s is equivalent to the study of such type of curvature restriction on the conharmonic curvature tensor $K$. The study of a commutative semisymmetric type and commutative pseudosymmetric type curvature restrictions on any $B$-tensor of class 2 is equivalent to the study of such type of restrictions on the conformal curvature tensor $C$. Moreover, each commutative and first type curvature restrictions on any $B$-tensor of class 2 with constant coefficients give rise only one structure i.e., the structure due to $K$. Also each commutative and second type curvature restrictions on any $B$-tensor of class 2 gives rise the same structure as to $C$.\\
(iii) study of a symmetric type and recurrent type curvature restrictions on any $B$-tensor of class 3 with constant $a_i$'s is equivalent to the study of such type of curvature restriction on the concircular curvature tensor $W$. Again the studies of locally symmetric, recurrent, commutative semisymmetric type and commutative pseudosymmetric type curvature restrictions on any $B$-tensor of class 3 are equivalent to the studies of such type of restrictions on the Riemann-Christoffel curvature tensor $R$. Moreover, each commutative and first type curvature restriction on any $B$-tensor of class 3 with constant coefficients give rise only one structure i.e., the structure due to $W$. Also each commutative and second type curvature restriction on any $B$-tensor of class 3 gives rise the same structure as to $R$.\\
(iv) study of a symmetric type and recurrent type curvature restrictions on any $B$-tensor of class 4 with constant $a_i$'s is equivalent to the study of such type of curvature restriction on the Riemann-Christoffel curvature tensor $R$. The study of a commutative semisymmetric type and commutative pseudosymmetric type curvature restrictions on any $B$-tensor of class 4 is equivalent to the study of such type of restrictions on the curvature tensor $R$. Moreover, each commutative and first type curvature restriction on any $B$-tensor of class 4 with constant coefficients give rise only one structure i.e., the structure due to $R$. Also each commutative and second type curvature restriction on any $B$-tensor of class 4 gives rise also the same structure as to $R$.\\
\indent Finally, we also conclude that for future study of any kind of curvature restriction (discussed earlier) on various curvature tensors, we have to study such curvature restriction on the tensor $B$ only and as a particular case we can obtain the results for various curvature tensors. We also mention that to study various curvature restrictions on the tensor $B$, we have to consider only the form of $B$ as given in (\ref{eqgen}) but not as the large form given in (\ref{tensor}).\\
\indent However, the problem of various structures for any two $B$-tensors of different class is still remain for further study.\\
%
\noindent \textbf{Acknowledgement:} The second named author gratefully acknowledges  to CSIR, New Delhi [File No. 09/025 (0194)/2010-EMR-I] for the financial assistance.
%



\begin{thebibliography}{99}\baselineskip=16pt
%
\bibitem{adm}
Adam$\grave{\mbox{o}}$w, A. and Deszcz, R., \emph{On totally umbilical submanifolds of some class of Riemannian manifolds}, Demonstratio Math., \textbf{16} (1983), 39-59.
\bibitem{CAH}
Cahen, M. and Parker, M., \emph{Sur des classes d'espaces pseudo-riemanniens symmetriques}, Bull. Soc. Math. Belg., {\bf 22} (1970), 339-354.
\bibitem{CAH1}
Cahen, M. and Parker, M., \emph{Pseudo-Riemannian symmetric spaces}, Mem. Amer. Math. Soc., {\bf 24} (1980).
\bibitem{ca}
Cartan, E., \emph{Sur une classe remarquable d'espaces de Riemannian}, Bull. Soc. Math. France, \textbf{54} (1926), 214- 264.
\bibitem{CHA}
Chaki, M. C., \emph{On pseudosymmetric manifolds}, An. Stiint. Ale Univ., AL. I. CUZA Din Iasi Romania, {\bf 33} (1987), 53-58.
\bibitem{LC}
Chongshan, L., \emph{On concircular transformations in Riemannian spaces}, J. Aust. Math. Soc. (Sr. A), {\bf 40} (1986), 218-225.
\bibitem{DEF}
Defever, F. and Deszcz, R., \emph{On semi-Riemannian manifolds satisfying the condition $R.R = Q(S,R)$}, Geometry and Topology of Submanifolds, III, World Sci., {\bf 1991}, 108-130.
\bibitem{DEF1}
Defever, F. and Deszcz, R., \emph{On warped product manifolds satisfying a certain curvature condition}, Atti. Acad. Peloritana Cl. Sci. Fis. Mat. Natur., {\bf 69} (1991), 213-236.
\bibitem{D0}
Defever, F., Deszcz, R., Hotlo$\acute{\mbox{s}}$, M., Kucharski, M. and Sent$\ddot{\mbox{u}}$rk, Z., \emph{Generalisations of Robertson-Walker spaces}, Annales Univ. Sci. Budapest. E$\ddot{\mbox{o}}$tv$\ddot{\mbox{o}}$s Sect. Math., {\bf 43} (2000), 13--24.
\bibitem{DA1}
Desai, P. and Amur, K., \emph{On $W$-recurrent spaces}, Tensor N. S., {\bf 29} (1975), 98-102.
\bibitem{DA2}
Desai, P. and Amur, K., \emph{On symmetric spaces}, Tensor N. S., {\bf 29} (1975), 185-199.
\bibitem{DR1}
Deszcz, R., \emph{Notes on totally umbilical submanifolds}, in: Geometry and Topology of Submanifolds, Luminy, May 1987, World Sci. Publ., Singapore, \textbf{1989}, 89-97.
\bibitem{DES}
Deszcz, R., \emph{On pseudosymmetric spaces}, Bull. Belg. Math. Soc., Series A, {\bf 44} (1992), 1-34.
\bibitem{D11}
Deszcz, R. and Glogowska, M., \emph{Some examples of nonsemisymmetric Ricci-semisymmetric 
hypersurfaces}, Colloq. Math., {\bf 94}(2002), 87--101.
\bibitem{D2}
Deszcz, R., Glogowska, M., Hotlo$\acute{\mbox{s}}$, M. and
Sent$\ddot{\mbox{u}}$rk, Z., \emph{On certain quasi-Einstein semi-symmetric hypersurfaces}, Annales Univ. Sci. Budapest. E$\ddot{\mbox{o}}$tv$\ddot{\mbox{o}}$s Sect. 
Math., {\bf 41}(1998), 151--164.
\bibitem{DGHS} Deszcz, R.,  G\l ogowska, M., Hotlo\'{s}, M. and Sawicz, K., 
\emph{A Survey on Generalized Einstein Metric Conditions},
Advances in Lorentzian Geometry:
Proceedings of the Lorentzian Geometry Conference in Berlin,
AMS/IP Studies in Advanced Mathematics {\bf{49}}, S.-T. Yau (series ed.),
M. Plaue, A.D. Rendall and M. Scherfner (eds.), 2011, 27--46.
\bibitem{des7}
Deszcz, R. and Grycak, W., \emph{On some class of warped product manifolds}, Bull. Inst. Math. Acad. Sinica, \textbf{15} (1987), 311-322.
\bibitem{D10}
Deszcz, R. and Hotlo$\acute{\mbox{s}}$, M., \emph{On hypersurfaces with type number two 
in spaces of constant curvature}, Annales Univ. Sci. 
Budapest. E$\ddot{\mbox{o}}$tv$\ddot{\mbox{o}}$s Sect. Math., {\bf
46} (2003), 19--34.
\bibitem{DUB}
Dubey, R. S. D., \emph{Generalized recurrent spaces}, Indian J. Pure Appl. Math., {\bf 10} (1979), 1508-1513.
\bibitem{sti}
Ewert-Krzemieniewski, S., \emph{On some generalisation of recurrent manifolds}, Math. Pannonica, \textbf{4/2} (1993), 191-203.
\bibitem{RG}
Garai, R. K., \emph{On recurrent spaces of first order}, Annali della Scuola Normale Superiore di Pisa - Classe di Scienze, \textbf{26(4)} (1972), 889-909.
\bibitem{Glodek}
Glodek, E., \emph{A note on riemannian spaces with recurrent projective curvature}, Pr. nauk. Inst. matem. i fiz. teor. Wr. Ser. stud. i mater. \textbf{1} (1970), 9-12.
\bibitem{GLOG}
Glogowska, M., \emph{Semi-Riemannian manifolds whose Weyl tensor is a Kulkarni-Nomizu square}, Publ. Inst. Math. (Beograd), {\bf 72:86}(2002), 95-106.
\bibitem{GLOG1}
Glogowska, M., \emph{On quasi-Einstein Cartan type hypersurfaces}, J. Geom. Phys., {\bf 58}(2008), 599-614.
\bibitem{BG}
Gupta, B., \emph{On projective-symmetric spaces}, J. Austral. Math. Soc., {\bf 4(1)}(1964), 113-121.
\bibitem{ishi}
Ishii, Y., \emph{On conharmonic transformation}, Tensor N. S., \textbf{7} (1957), 73-80.
\bibitem{Kow2} 
Kowalczyk, D., 
\emph{On the Reissner-Nordstr\"{o}m-de Sitter type spacetimes}, Tsukuba J. Math. {\bf{30}} (2006), 263--281.  
\bibitem{lee}
Lee, J. M., \emph{Riemannian Manifold, an Introduction to Curvature}, Springer (1997).
\bibitem{mik}
Mike\v{s}, J., \emph{Geodesic mappings of affine-connected and Riemannian spaces} (English), J. Math. Sci., New York \textbf{78(3)} (1996), 311-333.
\bibitem{mik2}
Mike\v{s}, J., Van\v{z}urov\'{a}, A. and Hinterleitner, I. \emph{Geodesic mappings and some generalisations}, p. 160, Olomouc \textbf{2009}.
\bibitem{ol} 
Olszak, K. and Olsak, Z., \emph{On pseudo-Riemannian manifolds with recurrent concircular curvature tensor}, Acta Math. Hungar., \textbf{137(1-2)} (2012), 64-71.
\bibitem{PV}
Petrovi\'c-Torgasev, M. and Verstraelen, L., \emph{On the concircular curvature tensor, the projective curvature tensor and the Einstein curvature tensor of Bochner-Kaehler manifolds}, Math. Rep. Toyama Univ., {\bf 10} (1987), 37-61.
\bibitem{pokh1}
Pokhariyal, G. P., \emph{Relativistic significance of Curvature tensors}, Int. J. Math. Math. Sci., \textbf{5(1)} (1982), 133-139.
\bibitem{pokh2}
Pokhariyal, G. P. and Mishra R. S., \emph{Curvature tensor and their relativistic significance}, Yokohama Math. J., \textbf{18} (1970), 105-108.
\bibitem{pokh3}
Pokhariyal, G. P. and Mishra R. S., \emph{Curvature tensor and their relativistic significance II}, Yokohama Math. J., \textbf{19(2)} (1971), 97-103.
\bibitem{prasad}
Prasad, B., \emph{A pseudo projective Curvature tensor on a Riemannian manifold}, Bull Calcutta Math. Soc., \textbf{94(3)} (2002), 163-166.
\bibitem{prasad2}
Prasad, B. and Maurya, A., \emph{Quasi concircular Curvature tensor on a Riemannian manifold}, Bull Calcutta Math. Soc., \textbf{30} (2007), 5-6.
\bibitem{prasad3}
Prasad, B., Doulo, K., and Pandey, P. N., \emph{Generalized quasi-conformal Curvature tensor on a Riemannian manifold}, Tensor N.S., \textbf{73} (2011), 188-197.
\bibitem{RL}
Rahaman, M. S. and Lal, S., \emph{On the concircular curvature tensor of Riemannian manifolds}, International Centre of for Theoretical Physics, \textbf{1990}.
\bibitem{RT}
Reynolds, R. F. and Thompson, A. H., \emph{Projective-symmetric spaces}, J. Australian Math. Soc., \textbf{7(1)} (1967), 48-54.
\bibitem{Ru1}
Ruse, H. S., \emph{On simply harmonic spaces}, J. London Math. Soc., {\bf 21} (1946), 243-247.
\bibitem{Ru2}
Ruse, H. S., \emph{On simply harmonic `kappa spaces' of four dimensions}, Proc. London Math. Soc., {\bf 50} (1949), 317-329.
\bibitem{Ru3}
Ruse, H. S., \emph{Three dimensional spaces of recurrent curvature},Proc. London Math. Soc., {\bf 50} (1949), 438-446.
\bibitem{sel}
Selberg, A., \emph{Harmonic analysis and discontinuous groups in weakly symmetric Riemannian spaces with applications
to Dirichlet series}, Indian J. Math., \textbf{20} (1956), 47-87.
\bibitem{SD}
Shaikh, A. A., Deszcz, D., Hotlo\'s, M., Je\l owicki and Kundu, H.,\emph{On pseudosymmetric manifolds}, Preprint.
\bibitem{SJ}
Shaikh, A. A. and Jana, S. K.,\emph{A Pseudo quasi-conformal curvature tensor on a Riemannian manifold}, South East Asian J. Math. Math. Sci., \textbf{4(1)} (2005), 15-20.
\bibitem{SP}
Shaikh, A. A. and Patra, A., \emph{On a generalized class of recurrent manifolds}, Archivum Mathematicum, {\bf 46} (2010), 39-46.
\bibitem{ROY1}
Shaikh, A. A. and Roy, I., \emph{On quasi generalized recurrent manifolds}, Math. Pannonica, {\bf 21(2)} (2010), 251-263.
\bibitem{ROY}
Shaikh, A. A. and Roy, I., \emph{On weakly generalized recurrent manifolds}, Annales Univ. Sci. Budapest. E$\ddot{\mbox{o}}$tv$\ddot{\mbox{o}}$s Sect. Math., {\bf 54} (2011), 35-45.
\bibitem{singh}
Singh, J. P., \emph{On m-Projective Recurrent Riemannian Manifold}, Int. J. Math. Analy., \textbf{6(24)} (2012), 1173 - 1178.
\bibitem{soos}
So\'os, G., \emph{\"Uber die Geod\"atischen Abbildungen Von Riemannschen R\"aumen auf Projectiv-symmetrische Riemannsche R\"aume}, Acta Acad. Sci. Hungarica, \textbf{9} (1958), 359-361.
\bibitem{sz}
Szab$\acute{\mbox{o}}$, Z. I., \emph{Structure theorems on Riemannian spaces satisfying $R(X,Y).R=0$} I, The local version, J. Diff. Geom., \textbf{17} (1982), 531-582.
\bibitem{tac}
Tachibana, S.,  \emph{A Theorem on Riemannian manifolds of positive curvature operator}, Proc. Japan Acad., \textbf{50} (1974), 301-302.
\bibitem{tb} 
Tam$\acute{\mbox{a}}$ssy, L. and Binh, T. Q., \emph{On weakly symmetric and weakly projective symmetric Riemannian manifolds}, Coll. Math. Soc. J. Bolyai, \textbf{50} (1989), 663-670.
\bibitem{tri} 
Tripathi, M. M. and Gupta, P., \emph{$\tau$-curvature tensor on a semi-Riemannian manifold}, J. Adv. Math. Stud., \textbf{4(1)} (2011), 117-129.
\bibitem{Ag} 
Walker, A. G., \emph{On Ruse's spaces of recurrent curvature}, Proc. London Math. Soc., \textbf{52} (1950), 36-64.
\bibitem{yano1}
Yano, K., \emph{Concircular geometry, I-IV}, Proc. Imp. Acad. Tokyo, \textbf{16} (1940), 195-200, 354-360, 442-448, 505-511.
\bibitem{yano4}
Yano, K. and Sawaki S., \emph{Riemannian manifolds admitting a conformal transformation group}, J. Diff. Geom., \textbf{2} (1968), 161-184.
\end{thebibliography}
\end{document}